\documentclass[12pt, draftclsnofoot, onecolumn, romanappendices]{IEEEtran}

\usepackage{bm}
\usepackage{graphicx}
\usepackage{color}
\usepackage{caption}
\usepackage{amsmath,amsthm,amsfonts,amssymb}
\usepackage{listings}
\usepackage{here}
\usepackage{cite}
\usepackage{hyperref}
\usepackage{bookmark}
\usepackage{footnote}
\makesavenoteenv{tabular}
\makesavenoteenv{table}

\title{Holonomic Gradient Method-Based \\CDF Evaluation for the Largest Eigenvalue of \\a Complex Noncentral Wishart Matrix}
\author{
Fadil H. Danufane\thanks{F.~H.~Danufane was with Kanazawa University, Japan. He is now with the Research Center for Electronics and Telecommunications, Indonesian Institute of Sciences, Bandung, Indonesia.}, 
Katsuyoshi Ohara\thanks{K.~Ohara is with the Faculty of Mathematics and Physics, Kanazawa University, Japan}, 
Nobuki Takayama\thanks{N.~Takayama is with the Department of Mathematics, Kobe University, Japan.}, 
Constantin Siriteanu\thanks{C.~Siriteanu is with the Graduate School of Information Science and Technology, Osaka University, Japan.}
}

\date{}

\usepackage{ragged2e}
\usepackage{booktabs}

\def\pd#1{\partial_{#1}}

\def\comment#1{}
\def\C{\mathbf{C}}
\def\qed{\hskip 15pt /\kern -2pt/ \bigbreak}

\def\tr{\operatorname{tr}\,}
\def\pr{\operatorname{Pr}\,}

\newtheorem{lemma}{Lemma}
\newtheorem{theorem}{Theorem}

\newtheorem{definition}{Definition}

\newtheorem{remark}{Remark}

\newtheorem*{objective*}{Objective}

\newcommand{\mb}{\mathbf}

\newcommand{\uh}{{\mb{h}}}

\newcommand{\uf}{{\mb{f}}}

\newcommand{\uhd}{{\mb{h}_{\text{d}}}}

\newcommand{\ug}{{\mb{g}}}

\newcommand{\uPhi}{{\pmb{\Phi}}}

\newcommand{\ulambda}{{\pmb{\lambda}}}

\newcommand{\ur}{{\mb{r}}}

\newcommand{\uQ}{{\mb{Q}}}

\newcommand{\uH}{{\mb{H}}}
\newcommand{\uHH}{{\mb{H}}^{\sf H}}

\newcommand{\uG}{{\mb{G}}}

\newcommand{\uHd}{{\uH_{\text{d}}}}

\newcommand{\uHdH}{{\uH_{\text{d}}^{\sf H}}}

\newcommand{\uHdn}{{\uH_{\text{d,n}}}}

\newcommand{\uHr}{{\uH_{\text{r}}}}

\newcommand{\uHrn}{{\uH_{\text{r,n}}}}

\newcommand{\NT}{{N_{\text{T}}}}

\newcommand{\NR}{{N_{\text{R}}}}

\newcommand{\uP}{{\mb{P}}}
\newcommand{\uR}{{\mb{R}}}

\newcommand{\uRT}{{\mathbf{R}_{\text{T}}}}

\newcommand{\uRTK}{{\mathbf{R}_{\text{T,}K}}}

\newcommand{\uRTKinv}{{\mathbf{R}_{\text{T,}K}^{-1}}}

\newcommand{\uS}{{\mb{S}}}
\newcommand{\uSd}{{\uS_{\text{d}}}}

\newcommand{\ux}{{\mb{x}}}

\newcommand{\un}{{\mb{n}}}

\newcommand{\uwT}{{\mb{w}}_{\text{T}}}

\newcommand{\mzero}{\mb{0}}

\newcommand{\mbI}{\mb{I}}

\def\hF{{\kern.15em{}_{1}\kern-.05emF_{1}}}

\newcommand{\hFoo}{{\kern.15em{}_{1}\kern-.05em {\sf F}_{1}}}

\newcommand{\hFooa}{{\kern.15em{}_{0}\kern-.05em {\sf F}_{1}}}

\newcommand{\hFooam}{{\kern.15em{}_{0}\kern-.05em {\sf \widetilde{F}}_{1}}}

\def\arraystretch{1.1}

\begin{document}
\maketitle

\begin{abstract}
The outage probability of maximal-ratio combining (MRC) for a multiple-input multiple-output (MIMO) wireless communications system under Rician fading is given by the cumulative distribution function (CDF) for the largest eigenvalue of a complex noncentral Wishart matrix. This CDF has previously been expressed as a determinant whose elements are integrals of a confluent hypergeometric function. For the determinant elements, conventional evaluation approaches, e.g., truncation of infinite series ensuing from the hypergeometric function or numerical integration, can be unreliable and slow even for moderate antenna numbers and Rician $ K $-factor values. Therefore, herein, we derive by hand and by computer algebra also differential equations that are then solved from initial conditions  computed by conventional approaches. This is the holonomic gradient method (HGM). Previous HGM-based evaluations of MIMO relied on differential equations that were not theoretically guaranteed to converge, and, thus, yielded reliable results only for few antennas or moderate $ K $. Herein, we reveal that gauge transformations can yield differential equations that are {\emph{stabile}}, i.e., guarantee HGM convergence. The ensuing HGM-based CDF evaluation is demonstrated reliable, accurate, and expeditious in computing the MRC outage probability even for very large antenna numbers and values of $ K $. 
\end{abstract}

\section{Introduction}

\subsection{Background and Model}

For multiple-input multiple-output (MIMO) wireless communications systems, maximal-ratio combining (MRC) is often considered for analysis and implementation\cite{Kang_jsac_03}\cite{wu_twc_16}. In MIMO MRC, each data symbol is transmitted into and received from the dominant singular mode of the $ \NR \times \NT $ channel matrix $ \uH $. This relatively low-complexity technique\footnote{Only knowledge of the dominant right and left singular vectors of $ \uH $ is needed, at the transmitter and receiver, respectively.} maximizes the symbol-detection signal-to-noise ratio (SNR) and, thus, the diversity gain. 

The MRC SNR is proportional to the largest eigenvalue $ \phi_s $ of $ \uHH \uH $\cite[Eq.~(27)]{Kang_jsac_03}. Therefore, the MRC outage probability, i.e., the probability of failing to achieve an SNR that ensures reliable symbol detection\cite{Kang_jsac_03}\cite{wu_twc_16}, is given by the cumulative distribution function (CDF) of $ \phi_s $. 

For complex circularly-symmetric Gaussian-distributed $ \uH $ with nonzero mean $ \uHd $, i.e., for Rician fading, and zero row correlation, matrix $ \uHH \uH $ has the complex noncentral Wishart distribution\cite{Kang_jsac_03}\cite{mckay_isit_05}. Then, the distribution of $ \phi_s $ has been characterized several times before, for various settings, with expressions that are difficult to evaluate reliably, as discussed next. 

\subsection{Previous Work and Its Limitations}

For the largest eigenvalue of a central Wishart matrix, Hashiguchi {\it{et al.}}\cite{hashiguchi_jma_2013} computed the CDF expressed as a confluent hypergeometric function $ {{}_1F_1}(\cdot; \cdot; \cdot) $ with matrix third argument\cite[Eq.~(4)]{hashiguchi_jma_2013}. 
Because here we are interested in the complex-valued noncentral Wishart distribution, we shall start from Kang and Alouini's paper\cite{Kang_jsac_03}. For the case when mean $ \uHd $ has arbitrary rank but zero row and column correlation, their\cite[Eq.~(2)]{Kang_jsac_03} is a determinantal expression for the CDF of the largest eigenvalue $ \phi_s $ of $ (K + 1) \uHH \uH $, where $ K $ is the Rician factor.
The CDF determinant elements are integrals of the special-case confluent hypergeometric function $ {{}_0F_1}(\cdot; \cdot; \lambda y) $, where $ \lambda $ stands for the eigenvalues of $ (K + 1) \uHdH \uHd $, and $ y $ is the integration variable. However, because for the eigenvalues $ \lambda_i $ of  $ (K + 1) \uHdH \uHd $ the typical normalization yields $ \sum \lambda_i = K \NR \NT $, larger $ K $, $ \NR $, or $ \NT $ implies the requirement to evaluate $ {{}_0F_1}(\cdot; \cdot; \lambda y) $ for larger $ \lambda $.
 
The obvious approach to compute the elements of the CDF determinant is numerical integration of native implementations of $ {{}_0F_1}(\cdot; \cdot; \lambda y) $ in software tools such as {Mathematica} or Matlab\footnote{Technical details of such implementations are not publicly available.}.
Another approach is truncation of the infinite series ensuing from the well-known formula\cite[Eq.~(13.2.2), p.~322]{NIST_book_10} for $ {{}_0F_1}(\cdot; \cdot; \lambda y) $. 
However, as $ \lambda $ becomes larger, these conventional approaches can lead to numerical inaccuracy (and even divergence) and/or long computation time, as illustrated in this paper. 
Note that\cite{Kang_jsac_03} showed MRC outage probability results  only for small antenna numbers, i.e., $ \max (\NR + \NT) = 5 $, and small Rician factor, i.e., $ K = 0 $~dB.

For $ \uH $ with arbitrary-rank $ \uHd $ and uncorrelated columns and rows, Wu {\it{et al.}} have recently expressed the CDF of the dominant eigenvalue $ \phi_s $ of $ \uHH \uH $ in\cite[Eq.~(3)]{wu_twc_16}. However, their  expression is an infinite series whose coefficients $ \alpha_k $, $ k = 0, 1, 2, \ldots $, expressed in\cite[Eq.~(4)]{wu_twc_16}, 
are difficult to express beyond $ k = 2 $. Further, the series was approximated with only its first term in\cite[Eq.~(19)]{wu_twc_16}, and MRC outage probability results were shown only for small antenna numbers, i.e., $ \max (\NR + \NT) = 6 $, and small-to-moderate Rician factor, i.e., $ K \leq 10 $~dB.

Furthermore, these MIMO performance evaluation limitations are not MRC-specific. As another example, for MIMO multiplexing with zero-forcing (ZF) detection, which requires knowledge of $ \uH $ only at the receiver, Siriteanu {\it{et al.}} analyzed and evaluated the performance for Rician fading with $ \text{rank}(\uHd) = 1 $ for a special case in\cite{Siriteanu_twc_14} and the general case in\cite{Siriteanu_twc_15}. 
The moment generating function of the ZF symbol-detection SNR, which for stream $ i = 1: \NT$ is proportional to $1/[(\uHH \uH)^{-1}]_{i,i}$, was characterized in\cite[Eq.~(19)]{Siriteanu_twc_14}\cite[Eq.~(52)]{Siriteanu_twc_15} with the confluent hypergeometric function $ {{}_1F_1}(\cdot; \cdot; \cdot) $ of third argument proportional with $ K \NR \NT $\cite[Eq.~(66)]{Siriteanu_twc_15}. The truncation of ensuing infinite series that characterize the ZF performance was shown to diverge numerically\footnote{This is because the infinite series for $ {{}_1F_1}(\cdot; \cdot; \cdot) $ diverges numerically for large values of the third argument\cite{muller_nm_01}.} for $ \NR = 6 $, $ \NT = 4 $, and moderate $ K $ values in\cite{Siriteanu_twc_14}\cite{Siriteanu_twc_15}. 

Such limitations are to be expected when evaluating any MIMO transceiver technique for Rician fading because the distribution of matrix $ \uHH \uH $, which determines the performance, is then \underline{noncentral} Wishart, i.e., characterized by a confluent hypergeometric function with matrix argument\cite[Eq.~(3)]{mckay_isit_05} that is inherently difficult to manipulate and compute\cite{koev_mc_06}.

\subsection{{Alternative Computation by the Holonomic Gradient Method} (HGM)}

As mentioned above, Siriteanu {\it{et al.}}\cite{Siriteanu_twc_14}\cite{Siriteanu_twc_15} evaluated MIMO ZF performance for Rician fading with $ \text{rank}(\uHd) = 1 $. First, they derived, by hand in\cite{Siriteanu_twc_14} and by computer algebra in\cite{Siriteanu_twc_15}, linear ordinary differential equations (LODEs) satisfied by infinite series that describe ZF performance measures, e.g., SNR probability density function in\cite{Siriteanu_twc_14} and outage probability in\cite{Siriteanu_twc_15}.
Then, they numerically evaluated these measures by a hybrid approach that solves LODEs for a desired value of $ K $ by starting from initial conditions computed by truncating infinite series for a value of $ K $ that is sufficiently small to ensure numerical convergence for the truncation. This approach has become known as the \emph{holonomic gradient method} (HGM)\cite{Siriteanu_twc_14}\cite{Siriteanu_twc_15}. 
However, HGM was applied in\cite{hashiguchi_jma_2013} to compute the largest-eigenvalue CDF only for a central Wishart matrix (i.e., for $ K = 0 $).
On the other hand, the HGM-based evaluations of MIMO ZF for Rician fading from\cite{Siriteanu_twc_14}\cite{Siriteanu_twc_15} were found reliable only for small MIMO deployments (e.g., $ \NR = 6 $, $ \NT = 4 $) or for small-to-moderate values of $ K $ ($<10$~dB).

HGM introductions, applications, and references appear in\cite{Siriteanu_twc_14}\cite{Siriteanu_twc_15}\cite[Ch.~6]{hibi_book_13}\cite{HGM_references}. HGM entails concepts, methods, and terminology from differential equations, abstract algebra, and Gr\"obner basis (or generator) computation, e.g., {\it{holonomic function and system}}, {\it{Pfaffian system}}, {\it{ring}}, {\it{ideal}}, etc., which are used herein but are thoroughly introduced in\cite[Ch.~1]{hibi_book_13}\cite{buchberger_mssp_01}\cite{sturmfels_nams_05}. Gr\"obner bases applications in signal processing are presented in\cite{lin_tcs_08}.

\subsection{Contribution: Reliable HGM-Based MIMO MRC Evaluation}
Herein, for the evaluation of the CDF of $ \phi_s $, i.e., the MRC outage probability, for uncorrelated Rician fading with $ \uHd $ of arbitrary rank, we propose a reliable HGM-based approach that avoids the limitations of numerical integration and series truncation, as well as limitations of the HGM-based MIMO ZF evaluations from\cite{Siriteanu_twc_14}\cite{Siriteanu_twc_15}. 
Our contributions are as follows:

\begin{enumerate}
\item From Kang and Alouini's integral expression of the CDF determinant elements\cite[Eq.~(2)]{Kang_jsac_03}, \underline{using the infinite series} known for $ {{}_0F_1}(\cdot; \cdot; \lambda y) $, we derive an infinite series and then deduce LODEs it satisfies from Gr\"obner bases obtained by computer algebra. Thereafter, we find that the ensuing HGM is reliable only for small $ \lambda $,  as for ZF in\cite{Siriteanu_twc_14}\cite{Siriteanu_twc_15}. We reveal the cause of this limitation to be that the used LODEs are not {\emph{stabile}}\footnote{Concept defined later.}.
\item From the integral expression, \underline{using the differential equation} known for $ {{}_0F_1}(\cdot; \cdot; \cdot) $, we derive a new LODE system and prove it {\emph{stabile}}, thus theoretically ensuring HGM effectiveness. This guarantees reliable HGM-based evaluation of MRC under Rician fading.
\item We reveal that the Step-2 stabile LODEs can also be obtained by a {\emph{gauge transformation}}\cite{jimbo_np_81} from the Step-1 LODEs, and that stabile LODEs can be obtained algorithmically by gauge transformations of other LODEs, {\underline{in general}}. 
\item We illustrate numerically that, unlike conventional integration and series truncation, HGM based on the Step-2 stabile LODEs, and further enhanced by gauge transformations, yields reliable evaluation of the CDF of $ \phi_s $ even for very large values of $ \lambda $, i.e., of the MRC outage probability for very large $ \NR $, $ \NT $, and Rician $ K $-factor. HGM is also much faster than integration and Monte-Carlo simulation.
\end{enumerate}

\subsection{Notation}
\label{section_notation}
\begin{itemize}
\item Scalars, vectors, and matrices are represented with lowercase italics, lowercase boldface, and uppercase boldface, respectively, e.g.,  $ \lambda $, $ \uh $, and $ \uH $; 
superscripts $ \cdot^{\sf T} $ and $ \cdot^{\sf H} $ stand for transpose and Hermitian (i.e., complex-conjugate) transpose; $ \mbI_N $ is the $ N \times N $ identity matrix.
\item $ i = 1 : N $ stands for the enumeration $ i = 1, \, 2, \, \ldots, \, N $. 
\item $ \propto $ stands for `proportional to'; $ \approx $ stands for `approximately equal'.
\item $ \uh \sim {\cal{CN}}_{\NR} \left( \uhd, \uR \right) $ denotes an $ \NR \times 1 $ complex-valued circularly-symmetric Gaussian vector with mean $ \uhd $ and covariance matrix $ \uR $; an $ \NR \times \NT $ complex-valued circularly-symmetric Gaussian random matrix with mean $ \uHd $, row covariance $ \mbI_{\NR} $, and column covariance $ \uRT $, i.e., a matrix whose vectorized form is distributed as $ \text{vec}(\uHH) \sim {\cal{CN}}_{\NR\NT} \left( \text{vec}(\uHdH),  \mbI_{\NR} \otimes \uRT \right) $, is denoted herein as $ \uH \sim {\cal{CN}}_{\NR,\NT} \left( \uHd, \mbI_{\NR} \otimes \uRT \right) $, based on the definition from\cite[Sec.~I.F]{Siriteanu_twc_15}; subscripts $ \cdot_{\text{d}} $ and $ \cdot_{\text{r}} $ identify, respectively, deterministic and random components; subscript $ \cdot_{\text{n}} $ indicates a normalized variable; $ \mathbb{E} \{ \cdot \} $ denotes statistical average; 
\item $\frac{\partial^{N}}{\partial x_i^{N}} f(x_1, \ldots, x_n) $ stands for the $ N $th partial derivative with respect to variable $ x_i $ of function $ f $; $\theta_{x} = x\frac{\pd{}}{\pd{}x} $ is the {\emph{theta}} differential operator w.r.t.~$x$; $ l(x_1, \ldots, x_n) \bullet f $ represents the application of differential operator $ l $ to function $ f $.
\item $ {{}_0F_1}(\cdot; \cdot; \cdot) $ is the \textit{confluent hypergeometric function}\cite[Ch.~13]{NIST_book_10}; $ (n)_i $ is the Pochhammer symbol, i.e.,  $ (n)_0 = 1 $ and $ (n)_i = n (n + 1) \ldots (n + i - 1) = (n + i - 1)!/(n - 1)!$, $ \forall i \geq 1 $, and $ i! = i \cdot (i - 1) \cdots 2 \cdot 1$.
\end{itemize}

\subsection{Paper Organization}
Section~\ref{section_MIMO_MRC} introduces the MIMO signal and channel model as well as the MRC method and its SNR.
Section~\ref{sec:determinant} gives the determinantal expression derived by Kang and Alouini in\cite{Kang_jsac_03} for the CDF of  $ \phi_s $, i.e., the MRC outage probability. For the determinant elements, which are integrals of $ {{}_0F_1}(\cdot; \cdot; \cdot) $, we derive an infinite series. Section~\ref{sec:local-hgm} derives LODEs satisfied by this infinite series, by hand and by computer algebra, and then implements HGM. 
Section~\ref{sec:global-hgm} derives by hand alternative LODEs from the integral expression of the determinant elements, and proves them stabile. For an enhanced HGM based on these LODEs and gauge transformations, Section~\ref{sec:numerical_results} shows results for various antenna numbers, sets of eigenvalues of $ \uHd $, and values of $ K $.

\section{MIMO Model, MRC SNR, Outage Probability}
\label{section_MIMO_MRC}

We consider a point-to-point MIMO wireless communications link. The transmitter is equipped with an $ \NT $-element antenna array that transmits the vector $  b \frac{E_{\text{b}}}{\NT} \uwT $, where $ b $ is the complex-valued transmitted symbol taken from a constellation with unit average energy, $ \frac{E_{\text{b}}}{\NT} $ is the energy consumed per transmitting antenna, i.e., $ E_{\text{b}} $ is the energy consumed per transmitted symbol, and $ \uwT $ is the unit-norm transmit combiner (beamformer).

Herein, we assume that the transmitted signal vector encounters a fading radio channel characterized by the matrix $ \uH $ distributed as $ \uH \sim {\cal{CN}}_{\NR,\NT} \left( \uHd, \mbI_{\NR} \otimes \uRTK \right) $, with transmit-side correlation matrix given by $ \uRTK = \frac{1}{K+1} \mbI_{\NT} $. Thus, we assume Rician fading that is uncorrelated both at the transmitter and the receiver.
With its deterministic and random components denoted as $\uHd$ and $\uHr$, respectively, we can write the channel matrix as
\begin{eqnarray}
\label{equation channel_matrix_components}
\uH =  \uHd + \uHr = \sqrt{\frac{K}{K+1}} \, \uHdn + \sqrt{\frac{1}{K+1}} \, \uHrn,
\end{eqnarray}
where $ \uHdn $ and $ \uHrn $ are normalized as usual, i.e., 
\begin{eqnarray}
\label{equation_Hdn_norm_assumption}
\label{equation_Hrn_norm_assumption}
\| \uHdn \|^2 = \mathbb{E} \{ \| \uHrn \|^2 \}  =  \NR \NT, \, \text{which implies }  \mathbb{E} \{ \| \uH \|^2 \} = \NR \NT.
\end{eqnarray}
Above, $ K $ is the Rician factor, i.e.,
\begin{eqnarray}
\label{equation K_definition}
K = \frac{\| \uHd \|^2}{\mathbb{E} \{ \| \uHr \|^2 \} } = \frac{ \frac{K}{K+1} \| \uHdn \|^2}{\frac{1}{K+1} \mathbb{E} \{ \| \uHrn \|^2 \} }.
\end{eqnarray}
Finally, the receiver front-end adds the complex additive noise vector $\un \sim {\cal{CN}}_{\NR} (\mzero, N_0 \, \mbI_{\NR})$. Thus, the $ \NR $-dimensional received signal vector  is
\begin{align}
\label{equation_received_signal}
\ur = \sqrt{\frac{{E}_{\text{b}}}{\NT}} \, \uH \uwT b + \un.
\end{align}
We define the per-symbol transmit SNR as
\begin{align}
\label{equation_Gammas}
\Gamma_{\text{b}} = \frac{{E}_{\text{b}}}{N_0} \frac{1}{\NT}.
\end{align}

\begin{remark}
\label{lemma_MRC_SNR}
For the MIMO signal model in~(\ref{equation_received_signal}), MRC maximizes the symbol-detection SNR by beamforming with the dominant right and left singular vectors of $ \uH $ at the transmitter and receiver, respectively. 
Consequently, the MRC SNR is given by $ \Gamma_{\text{b}} $ multiplied with the maximum eigenvalue of matrix $ \uHH \uH $ --- see\cite[Sec.~III.B]{Kang_jsac_03} for the simple proof.
\end{remark}

In general, the outage probability is defined as the probability that the symbol-detection error probability --- instantaneous, i.e., given $ \uH $ --- exceeds a certain threshold, i.e., $ P_\text{e} \ge P_\text{e,th} $. 
The threshold SNR that yields the threshold error probability is denoted herein with $ \Gamma_{\text{th}} $. For instance, for QPSK modulation, a symbol error probability level of $ 10^{-2} $ can be achieved for $ \Gamma_{\text{th}} \approx 8.2 $~dB. 
Thus, the outage probability of a detection method can also be defined as the probability that its symbol-detection SNR does not exceed the threshold value $ \Gamma_{\text{th}} $. Clearly, the outage probability depends on the distribution of symbol-detection SNR, which is determined by the distribution of the MIMO channel matrix $ \uH $.

In the next section, we shall characterize the MRC outage probability based on Kang and Alouini's work from\cite{Kang_jsac_03}. 
They considered the matrix $ \uS = \uRTKinv \uHH \uH $, with eigenvalues\footnote{All these eigenvalues are real-valued and positive because, with probability 1, matrix $ \uS $ is Hermitian and has rank $s$.}  $\phi_1<\phi_2<\cdots<\phi_s$, and the corresponding matrix $ \uSd = \uRTKinv \uHdH \uHd $ with eigenvalues $\lambda_1<\lambda_2<\cdots<\lambda_s$, where $s=\min(\NT,\NR)$ and $t=\max(\NT,\NR)$. 
They characterized the CDF of the dominant eigenvalue $ \phi_s $ of $ \uS $ based on the eigenvalues of $ \uSd $.
Nevertheless, they later assumed zero transmit-side correlation in order to express the MRC SNR and its outage probability.

\begin{remark}
\label{lemma_MRC_SNR_cor}
The MRC SNR characterized in Remark~\ref{lemma_MRC_SNR}, is given for $ \uRTK = \frac{1}{K+1} \mbI_{\NT} $ by the dominant eigenvalue $ \phi_s$ of $ \uS = \uRTKinv \uHH \uH = (K+1) \uHH \uH $ as follows 
\begin{eqnarray}
\label{equation_MRC_rho}
\rho_{\text{MRC}} = \frac{\Gamma_{\text{b}}}{K+1} \phi_s,
\end{eqnarray}
and the MRC outage probability is 
\begin{eqnarray}
\label{equation_MRC_Outage_Probability}
\pr(\rho_{\text{MRC}} \leq \Gamma_{\text{th}}) = \pr \left(\phi_s \leq \frac{(K+1) \Gamma_{\text{th}}}{\Gamma_{\text{b}}} \right).
\end{eqnarray}
\end{remark}

Expression~(\ref{equation_MRC_Outage_Probability}) confirms that a larger per-symbol transmit SNR $ \Gamma_{\text{b}} $ (i.e., a larger transmitted energy) yields a lower outage probability whereas a larger threshold SNR $ \Gamma_{\text{th}}  $ (i.e., a higher quality of service) yields a higher outage probability.
It also confirms that the MRC outage probability depends on the distribution of $ \phi_s $, which is characterized next, based on\cite{Kang_jsac_03}.

\section{CDF of Largest Eigenvalue of Complex Noncentral Wishart Matrix}
\label{sec:determinant}

We begin with results from\cite{Kang_jsac_03} that apply for the distribution $ \uH \sim {\cal{CN}}_{\NR,\NT} \left( \uHd, \mbI_{\NR} \otimes \uRTK \right) $. Then, for analysis tractability, we specialize to the case with $ \uH \sim {\cal{CN}}_{\NR,\NT} \left( \uHd, \mbI_{\NR} \otimes \frac{1}{K+1} \mbI_{\NT} \right) $.

\begin{definition}
\label{theorem_Wishart_distribution}
If $ \uH \sim {\cal{CN}}_{\NR,\NT} \left( \uHd, \mbI_{\NR} \otimes \uRTK \right) $, then the $ \NT \times \NT $ matrix $ \uHH \uH $ has a complex noncentral Wishart distribution with the following parameters: dimension $ \NT$, degrees of freedom $ \NR $, and noncentrality matrix $ \uRTKinv \uHdH \uHd  $\cite{Kang_jsac_03}\cite{mckay_isit_05}\cite{Siriteanu_twc_15}.
\end{definition}

Now we state the result proved by Kang and Alouini in\cite{Kang_jsac_03}.

\begin{theorem}
\label{theorem_Kang_Alouini_result}
If $ \uH \sim {\cal{CN}}_{\NR,\NT} \left( \uHd, \mbI_{\NR} \otimes \uRTK \right) $ and the $ \NT \times \NT $ matrix $ \uSd = \uRTKinv \uHdH \uHd   $ has $s$ nonzero distinct eigenvalues $0<\lambda_1<\lambda_2<\cdots<\lambda_s$, then the CDF of the largest eigenvalue $\phi_s$ of $ \uS =  \uRTKinv \uHH \uH $ evaluated at certain threshold value
$x$ is given by the determinantal expression
\begin{eqnarray}
\label{equation_CDF_det_original}
\pr(\phi_s \leq x)=\frac{e^{-(\lambda_1 + \cdots + \lambda_s)} }{ \left[\prod_{1\le i<j \le s} (\lambda_i - \lambda_j) \right] \left[ (t-s)!\right]^{s}} \det \uPhi(x),
\end{eqnarray}
where $\uPhi(x)$ is an $ s \times s $ matrix whose $(i,j)$th element is given by
\begin{eqnarray}\label{eq:integral-rep}
 [\uPhi(x)]_{i,j} =\int_0^x y^{t-i}e^{-y}{}_0F_1(;t-s+1; \lambda_j y) \mathrm{d} y.
\end{eqnarray}
\end{theorem}    

Recall from Remark~\ref{lemma_MRC_SNR} that the MRC SNR is given by the dominant eigenvalue of $ \uHH \uH $, regardless of the column-correlation matrix $ \uRTK $. However, because the dependence between the eigenvalues of $ \uHH \uH $ and  $ \uS =  \uRTKinv \uHH \uH $ is not straightforward for $ \uRTK \not \propto \mbI_{\NT} $, herein, we employ Theorem~\ref{theorem_Kang_Alouini_result} for MRC analysis only for\footnote{As in\cite{Kang_jsac_03}, where Kang and Alouini first derived the CDF for the dominant eigenvalue of $ \uS =  \uRTKinv \uHH \uH $ in\cite[Eq.~(2)]{Kang_jsac_03}, and then assumed zero correlation between the columns of $ \uH $ in order to express the MRC outage probability in\cite[Eq.~(29)]{Kang_jsac_03}.} $ \uRTK = \frac{1}{K+1} \mbI_{\NT} $. 
Then, Theorem~\ref{theorem_Kang_Alouini_result} characterizes the CDF of the dominant eigenvalue $ \phi_s $ which enters the MRC SNR as in~(\ref{equation_MRC_rho}), and determines the MRC outage probability as in~(\ref{equation_MRC_Outage_Probability}).

\begin{remark}
\label{remark_trace_KNTNR}
\label{remark_trace_KNTNR_1}
Our model and assumptions yield $\tr(\uSd) = \| (K + 1) \uHd \|^2 = K \NR \NT $, i.e., the eigenvalues $ \lambda_i $, $ i = 1 : s $, of $ \uSd $ satisfy $ \sum_{i=1}^{s} \lambda_i = K \NR \NT $.
Consequently, even a relatively small MIMO deployment with $ (\NT, \NR) = (4, 5) $, in moderate Rician fading with $ K = 7 $~dB (i.e., $ K = 5 $) yields $ \sum_{i = 1}^{4} \lambda_i  = 100 $. Much larger values for $ \lambda_i $'s ensue with increasing $ \NT, \NR $, and $ K $.  For instance, a truly massive MIMO system with $ (\NT, \NR) = (1000, 1000) $ in strong Rician fading with $ K = 20 $~dB, yields $ \sum_{i = 1}^{1000} \lambda_i = 10^8 $. 
\end{remark}

\begin{objective*}
Given the eigenvalues $ \lambda_i $, $ i = 1 : s $, and a threshold value $x$, we seek the ability to compute the value of the CDF $\pr(\phi_s\leq{}x)$ by first evaluating the function in~(\ref{eq:integral-rep}) and then the determinant in~(\ref{equation_CDF_det_original}). Finally, scaling as in~(\ref{equation_MRC_Outage_Probability}) yields the MRC outage probability.
\end{objective*}

Defining $k=t-i$ and $n=t-s+1$, the determinant entry from~(\ref{eq:integral-rep}) can be rewritten as
\begin{eqnarray}
\label{equation_Hnk_definition_integral}
H^k_n(x,\lambda) =\int_0^x y^{k}e^{-y} {}_0F_1(;n; \lambda y)\mathrm{d} y.
\end{eqnarray}
One may then attempt to compute $ H^k_n(x,\lambda)  $ by numerical integration, e.g., in Mathematica.
However, we show later that this approach is not time-efficient for realistically-large values of $ \lambda $.
Therefore, next, we derive an alternative $ H^k_n(x,\lambda) $ expression and discuss its computation.

By employing for the confluent hypergeometric function $ {}_0F_1(;n; \lambda_j y)  $ its well-known infinite series $ \sum_{i=0}^{\infty}\frac{1}{(n)_i}\frac{{(\lambda_j y)^i}}{{i}!}$\cite[Eq.~(13.2.2), p.~322]{NIST_book_10} and then changing integration and summation order  in~(\ref{equation_Hnk_definition_integral}) we have obtained the following infinite series 
\begin{equation}\label{eq:series-rep}
 H^k_n(x,\lambda)  = \frac{x^{k+1}}{k+1} \sum_{p=0}^\infty \sum_{q=0}^\infty \frac{ (-1)^q (k+1)_{p+q}}{(n)_p(1)_p(1)_q(k+2)_{p+q}}x^{p+q} \lambda^p,
\end{equation}
which can be proven to converge $ \forall  (x,\lambda) $.
One may attempt to compute $H^{k}_{n} (x,\lambda) $ with truncation 
\begin{equation}\label{eq:truncated}
 H^k_n(x,\lambda)  \approx H^k_{n,N} (x,\lambda)  = \frac{x^{k+1}}{k+1} \sum_{p=0}^N \sum_{q=0}^N \frac{ (-1)^q (k+1)_{p+q}}{(n)_p(1)_p(1)_q(k+2)_{p+q}}x^{p+q} \lambda^p,
\end{equation}
with $N$ determined as follows: given the pair $ (x,\lambda) $ and a positive scalar $\varepsilon$, we select $N$ as the smallest integer that yields numerical convergence, i.e., satisfies
$|H^k_{n,N+1} (x,\lambda)  - H^k_{n,N} (x,\lambda) | < \varepsilon\,| H^k_{n,N} (x,\lambda) |$.

Thus, we have attempted to compute $ H^k_n(x,\lambda) $ by implementing the series truncation~(\ref{eq:truncated}) in C with {\emph{double-precision arithmetic}}\footnote{The series truncation computations presented here have been executed on the machine described in Section~\ref{sec:numerical_results} on page~\pageref{computer_description}.}, and compared for correctness and execution time against computations in Risa/Asir with rational arithmetic\cite{Risa_software_17}. 
Thus, we have found that the double-arithmetic C implementation of series truncation~(\ref{eq:truncated}), with $ \varepsilon = 10^{-10} $, can compute accurately $ H^2_3(x,\lambda) $ for $ x = \lambda $ taking increasing values in the interval $ [1.2, 10] $; numerical convergence within $ \varepsilon = 10^{-10} $ is then achieved for $ N $ values that increase in the interval $ [14, 48] $, but the computation takes negligible time (i.e., less than a microsecond). 

However, for larger values of $ x $ and $ \lambda $ (which may arise in practice, cf.~Remark~\ref{remark_trace_KNTNR}), the double-arithmetic C implementation of series truncation~(\ref{eq:truncated}) can experience computational difficulties. For example, we have found that the computation of $ H^2_3(30,30) $ begins to break down for $ N = 85 $ because denominators in the terms of series~(\ref{eq:truncated}) reach the upper limit of double-precision representation. Consequently, for $ N > 85 $, these denominators are treated as `infinite' and subsequent series terms remain constant, which precludes numerical convergence. The double-arithmetic C computation of series truncation~(\ref{eq:truncated}) with $ N = 85 $ takes about 45 microseconds (but does not achieve numerical convergence).
On the other hand, the rational-arithmetic Risa/Asir implementation of series truncation~(\ref{eq:truncated}) successfully achieves numerical convergence within $ \varepsilon = 10^{-10} $ for $ N = 120 $, but takes $ 556 $ milliseconds. 

These limitations of integration and truncation have motivated us to seek more reliable and efficient alternative approaches.
Thus, we propose below new $ H^k_n(x,\lambda) $ evaluation approaches based on differential equations and HGM. 
First, in Section~\ref{sec:local-hgm}, we derive a $4$-dimensional system of LODEs for $ H^k_n(x,\lambda) $ from the infinite series in~(\ref{eq:series-rep}) based on Gr\"obner basis computation, but find the ensuing HGM ineffective.
Then, in Section~\ref{sec:global-hgm_all}, we use the integral relationship~(\ref{equation_Hnk_definition_integral}) to derive a $3$-dimensional system of LODEs for $  H^k_n(x,\lambda) $ and propose a series of HGM enhancements that lead to accurate CDF evaluation.

\section{HGM for $H^k_n(x,\lambda)$ Using Differential Equations for Infinite Series~(\ref{eq:series-rep})}
\label{sec:local-hgm}

\subsection{Holonomic Functions and the Holonomic Gradient Method (HGM)}
    
HGM is a relatively new approach to numerically evaluate infinite series, such as our series for $H^k_n(x,\lambda)$ from~(\ref{eq:series-rep}). It relies on the notion of holonomic function and on the computation of Gr\"obner bases. 

\begin{definition}
\label{definition_holonomic}
A complex-valued function $f(\ux)$, where $ \ux=(x_1, \ldots, x_n)\in \mathbf{C}^n$, is {\it holonomic} if and only if,  $ \forall i = 1 : n $, there exists a LODE of some order $N_i$ w.r.t.~$ x_i $, as follows
\begin{eqnarray}
\bigg\{p_{i,N_i}(\ux) \frac{\partial^{N_i}}{\partial x_i^{N_i}}+\cdots+p_{i,1}(\ux)\frac{\partial}{\partial x_i}+p_{i,0}(\ux) \bigg\}
\bullet f(\ux) = 0,
\end{eqnarray}
where $p_{i,N_i}(\ux),\ldots,p_{i,0}(\ux)$ are polynomials.
The set of LODEs w.r.t.~all variables is called the holonomic system satisfied by $ f(x) $.
\end{definition}

\begin{definition}
\label{definition_Pfaffian}
A Pfaffian system is a set of LODEs of the form
\begin{eqnarray}
\label{equation_Pfaffian_ODES}
 \frac{\partial}{\partial x_i}  {\uf}(\ux)= \uP_i(\ux) {\uf}(\ux), \quad i=1, \ldots, n,
 \end{eqnarray}
where ${\uf}$ is a function vector, and matrices $ \uP_i $ 
satisfy the compatibility (integrability) condition $  \frac{\partial}{\partial x_j}  \uP_i =  \frac{\partial}{\partial x_i}  \uP_j $, $ \forall i, j = 1 : n. $\cite[p.~296]{hibi_book_13}. 
\end{definition}

Pfaffian LODE systems such as~(\ref{equation_Pfaffian_ODES}) are desirable because they can be solved numerically, e.g., by Runge--Kutta method, given the initial condition value $ {\uf}(\ux_0) $ for some $ \ux_0 $. 

\begin{definition}
\label{definition_HGM}
HGM for the computation of a function $ {\uf}(\ux) $ entails the following steps\cite[Ch.~6]{hibi_book_13}: 
\begin{enumerate}
\item Derive a holonomic system for $ {\uf}(\ux) $, e.g., by hand or by computer algebra.
\item Recast the holonomic system as a Pfaffian system.
\item Determine an initial condition $ {\uf}(\ux_0) $ for some $ \ux_0 $, e.g., by series truncation.
\item Numerically evaluate $ {\uf}(\ux) $ at desired point $ \ux $, e.g., by the Runge--Kutta method.
\end{enumerate}
\end{definition}

\subsection{Differential Equations and Pfaffian System from Infinite Series~(\ref{eq:series-rep}) for $H^k_n(x,\lambda)$}

\begin{theorem} 
\label{theorem_LODEs_4D}
\label{theorem_4d_LODE_pfaffians}
Function $H^k_n(x,\lambda)$ represented by the infinite series in~(\ref{eq:series-rep}) is holonomic because it satisfies the following LODEs  w.r.t.~$ x $ and~$ \lambda $, respectively,
\begin{eqnarray}
\label{eq:annihilating-x}
&& \big \{ \left[ (\theta_x+x-k-1)(\theta_x+x-k+n-2)-x\lambda \right]\theta_x \big \} \bullet  H^k_n(x,\lambda)  = 0,  \\
&& \big\{ \theta_{\lambda}^4+(\lambda+2n-4)\theta_{\lambda}^3+\left[ n^2-5n+5-(x+k+n) \lambda \right]\,\theta_{\lambda}^2  +\left[ x\lambda^2-(n-1)x \lambda \right. \nonumber \\
\label{eq:annihilating-lambda}
&& \quad \quad \quad \quad \left. -(k+1)(n-1) \lambda -(n-2)(n-1)\right] \theta_{\lambda} +(k \lambda+1)x \lambda \big \} \bullet  H^k_n(x,\lambda)  = 0, 
\end{eqnarray}
where $ \theta_{x} = x\frac{\pd{}}{\pd{}x} $ and $ \theta_{\lambda} =  \lambda \frac{\pd{}}{\pd{} \lambda } $.
These LODEs can be recast as the Pfaffian system
\begin{eqnarray}
\label{eq:pfaffian}
\frac{\pd{}}{\pd{}x}  \uf (x,\lambda)  & = & \uP (x,\lambda)  \, \uf (x,\lambda) , \quad \\
\label{eq:pfaffian_wrt_lambda}
\frac{\pd{}}{ \pd{}\lambda }  \uf (x,\lambda)  & = & \uQ (x,\lambda)  \,  \uf (x,\lambda) ,
\end{eqnarray}
where $ \uf $ is the $ 4 $-dimensional vector
\begin{eqnarray} 
\label{equation_vector_f_4D}
\uf  (x,\lambda) = \begin{matrix} (1 & \theta_{\lambda} & \theta_{\lambda}^2 & \theta_{\lambda}^3)^{\sf{T}} \end{matrix} \bullet u(x,\lambda),
\end{eqnarray}
with $ u (x,\lambda)  $ being a generic function\footnote{$ u(x,\lambda) $ stands for all four linearly independent solutions of the $ 4 $-dimensional system of LODEs in~(\ref{eq:pfaffian}).}, and
\begin{eqnarray}
\label{equation_P1_P2_Pfaffian}
\uP (x,\lambda) 
=
\frac{1}{x \lambda}
\left( 
\begin{array}{c@{\ \ }c@{\ \ }c@{\ \ }c}
a_1 & a_2 & -1 & 0 \\ 
0 & a_3 & a_2+1 & -1  \\ 
a_5 & a_6  & a_7 & n-1 \\ 
a_8 & a_9  & a_{10} & a_{11}
\end{array}
\right),
\,
\uQ  (x,\lambda) =\frac{1}{\lambda} \left( 
        \begin{array}{c@{\ }c@{\ }c@{\ }c}
            0 & 1 & 0 & 0\\
            0 & 0 & 1 & 0\\
            0 & 0 & 0 & 1\\
            -b_0 & -b_1 & -b_2& -b_3\\
        \end{array}
        \right),
\end{eqnarray}
where  $ a_i $ and $ b_j $ are polynomials in $ x $ and  $ \lambda $ shown in Table~\ref{equation_ai_bj}.
\renewcommand{\arraystretch}{0.8}
\begin{table}
\begin{center}
\caption{Elements of matrices $ \uP(x,\lambda) $ and $ \uQ (x,\lambda) $ from~(\ref{equation_P1_P2_Pfaffian})}
\label{equation_ai_bj}
\begin{tabular}{c | l}
           &  Expression \\ \hline
 $ a_1 $ &  $ (k+1) \lambda  $  \\
 $ a_2 $ &  $  \lambda -n+1 $  \\
 $ a_3 $ &  $ (k+1) \lambda +n-1 $  \\
 $ a_5 $ &  $ (k+1)x \lambda ^2 $  \\
 $ a_6 $ &  $ x \lambda ^2 -(n-1) \left[ x \lambda +(k+1) \lambda +n-1\right] $  \\
 $ a_7 $ &  $ -(x+n-1) \lambda +(n-1)(n-2)  $ \\
 $ a_8 $ &  $ -(n-2)(k+1)x \lambda ^2  $ \\
 $ a_9 $ &  $ (k-n+3)x \lambda ^2+(n-1)^2 \left[ x \lambda +(k+1) \lambda +n-1 \right]  $ \\
 $ a_{10} $ &  $ x \lambda ^2+(n-1)^2( \lambda -n+2)  $ \\
 $ a_{11} $ &  $ -x \lambda -(n-1)^2  $ \\ \hline
 $ b_3 $ &  $  \lambda +2n-4 $  \\
 $ b_2 $ &  $ n^2-5n+5-(x+k+n) \lambda  $ \\
 $ b_1 $ &  $ x \lambda ^2-(n-1)x \lambda -(k+1)(n-1) \lambda  -(n-2)(n-1)  $ \\
 $ b_0 $ &  $ (k \lambda +1)x \lambda  $ \\ \hline
\end{tabular}
\end{center}
\end{table}
\end{theorem}

\begin{IEEEproof}
Shown in Appendix~\ref{appendix_proof_theorem_ODE4}; it entails 1) by-hand deduction of two partial differential equations, and 2) LODEs deduction by computer-algebra-based computations of Gr\"obner bases.
\end{IEEEproof}

\subsection{Performance of HGM w.r.t.~$ \lambda $ Based on LODE System~(\ref{eq:pfaffian_wrt_lambda})}
\label{experiment_HGM_wrt_lambda}

Recall from Remark~\ref{remark_trace_KNTNR} that large values for $ \NT, \NR $, and $ K $ can yield large magnitudes for the eigenvalues $ \lambda_i $, $ i = 1 : s $, of $ \uSd $.
Thus, we have attempted to evaluate the function $H^k_n(x, \lambda )$ over a range of values for $ \lambda $ with HGM w.r.t.~$ \lambda $, based on system~(\ref{eq:pfaffian_wrt_lambda}). The numerical procedures and settings are summarized on the second line in Table~\ref{table_numerical_procedure} and relevant footnotes on page~\pageref{table_numerical_procedure}. 

Our HGM w.r.t.~$ \lambda $ based on system~(\ref{eq:pfaffian_wrt_lambda}) for $ \uf  (x,\lambda) $ comprises the following steps:
\begin{itemize}
\item Compute $  \uf  (x,\lambda) $ at some initial value $ \lambda_0 $ (i.e., the initial condition), by truncation~(\ref{eq:truncated}).
\item Solve the LODE system in~(\ref{eq:pfaffian_wrt_lambda}), and thus advance $  \uf  (x,\lambda) $ from initial value  $ \lambda_0 $ to the value of interest $ \lambda $, by the Runge--Kutta method.
\end{itemize}

We have attempted to compute $ H^2_3(1,\lambda) $ by the above HGM procedure from the initial condition $ H^2_3(1,\lambda_0 = 10^{-5}) $ obtained by series truncation, and compared with results obtained (upon numerical convergence) from the double-arithmetic C implementation of series truncation~(\ref{eq:truncated}). We have found that they agree closely when $ \lambda $ is near the initial value $ \lambda_0 $, but they begin to disagree as $ \lambda $ becomes larger, e.g., $ \lambda >10 $.  We conclude that the application of HGM based on the $ 4 $-dimensional LODE system from~(\ref{eq:pfaffian_wrt_lambda}) is reliable only locally, i.e., within a narrow range around small $ \lambda_0 $. The reason is explained further below.

\section{HGM for $H^k_n(x,\lambda)$ Using Stabile LODEs}
\label{sec:global-hgm_all}

\label{section_modified_HGM}
\label{sec:global-hgm}

We have seen above that HGM based on~(\ref{eq:pfaffian_wrt_lambda}) is reliable only locally.
This phenomenon can be understood based on the notion, defined below, of a LODE (system) that is {\emph{stabile}}\footnote{Which is different in meaning and spelling than the conventional notion of (Lyapunov) stable ODE.} for its solution.
In this section, we first directly derive a new LODE system stabile for $H^k_n(x,\lambda)$ from integration~(\ref{equation_Hnk_definition_integral}), and then propose {\emph{gauge transformations}} as a general means to obtain from a given LODE system that is not stabile a lower-dimensional LODE system that is stabile. 

\subsection{LODE Stabile for Its Solution}

\begin{definition}
Consider a holonomic function\footnote{Here, the generic variable symbol $ x $ stands for either of the variables of $ H^k_n(x,\lambda) $.} $ f(x) $ that satisfies a LODE of order $ m $, and denote with $ f_1(x), \ldots, f_m(x) $ its linearly independent solutions.
Then, let $f_i(x)$ be the dominant solution for $x \rightarrow \infty$, i.e., $|f_i(x)| \geq |f_j(x)|$, $ \forall j$.
We refer to the LODE as {\emph{stabile}}
for $f(x)$ if $ \lim_{x \rightarrow \infty} \frac{| f_i (x) |}{| f (x) |} < \infty $. (Note that a LODE is stabile or not regardless of the selected set of linearly independent solutions.) The notion of stabile LODE is defined analogously in the case of a vector-valued function,
by replacing $|\cdot|$ with a vector norm $|| \cdot ||$.
\end{definition}

For example, the solution space of LODE
$[{\mathrm{d}}^2/{\mathrm{d}}x^2 - 3 {\mathrm{d}}/{\mathrm{d}}x + 2] \bullet u = 0$
is spanned by $ f_1(x) = e^{2x} $, which is the dominant solution, and by $ f_2(x) = e^{x} $.
Consider the general LODE solution $f(x) = C_1 e^{2x} + C_2 e^{x} $. Then, if $C_1 \not= 0$ the LODE is stabile for $f$; otherwise, the LODE is not stabile for $f$.

\begin{remark}
\label{remark_stabile_not}
Typically, the initial condition $f(x_0)$ employed in evaluating function $f(x)$ by HGM is affected by some error (e.g., due to series truncation).
Then, HGM using a stabile LODE can theoretically\footnote{I.e., for computation on a machine with unlimited computational power and arithmetic precision.} still yield the function $ f(x) $ accurately as $ x $ increases, whereas HGM using a LODE that is not stabile may yield a  solution of much larger magnitude than $f(x)$. Thus, the latter can be effective in evaluating $f(x)$ only locally, for $ x $ around some small $ x_0 $.
\end{remark}

\begin{theorem} 
\label{th:not stabile-lambda-2017-12-05}
The LODE~(\ref{eq:annihilating-lambda}) and ensuing $4$-dimensional LODE system~(\ref{eq:pfaffian_wrt_lambda}) are not stabile for $H^k_n(x,\lambda)$
when $\lambda \rightarrow \infty$.
On the other hand, LODE~(\ref{eq:annihilating-x}) is stabile
whereas ensuing $4$-dimensional LODE system~(\ref{eq:pfaffian}) is not stabile when
$x \rightarrow \infty$.
\end{theorem}

\begin{IEEEproof}
See Appendix~\ref{section_proof_theorem_not stabile-lambda-2017-12-05}.
\end{IEEEproof}

\begin{remark}
The fact that LODE~(\ref{eq:annihilating-lambda}) and system~(\ref{eq:pfaffian_wrt_lambda}) are not stabile explains the ineffectiveness of the HGM w.r.t.~$ \lambda $ from Section~\ref{sec:local-hgm}. 
\end{remark}

\subsection{Direct Derivation From Integration~(\ref{equation_Hnk_definition_integral}) of New LODE System Stabile for $H^k_n(x,\lambda)$}

By differentiating~(\ref{equation_Hnk_definition_integral}) w.r.t.~$ x $, we obtain the following inhomogeneous equation for $ H^k_n(x,\lambda)  $:
    \begin{eqnarray}
    \label{eq:dxHkn_1}
        \frac{\pd{}}{\pd{}x} H^k_n(x,\lambda) = x^k e^{-x}{}_0F_1(;n;x \lambda).
    \end{eqnarray}
Now, for any integer $n>0$, there exists constant $c_n \neq 0 $ such that\cite{Dixit_book_15}
\begin{eqnarray}
\label{equation_0F1_approximation}
 \lim_{ z \rightarrow \infty} \frac{| {}_0F_1(;n; z) - c_n \, e^{2 \sqrt{z}}  (\sqrt{z})^{-n+1/2}|}{e^{2 \sqrt{z}}  (\sqrt{z})^{-n+1/2}} \rightarrow 0.
 \end{eqnarray}
Thus, for large  $ \lambda $, 
$ {}_0F_1(;n; x \lambda)$ increases with  $ \lambda $ as rapidly as $ e^{2 \sqrt{x \lambda}}  (\sqrt{x\lambda})^{-n+1/2}  $. 
We can avoid the exponential growth by 
working instead with the function (in which we assume $ \lambda $ to be a given constant)
\begin{eqnarray}
\label{equation_function_v}
v(x) =  e^{-2\sqrt {x \lambda} } {}_0F_1(;n; x \lambda)  =  e^x e^{-2\sqrt {x \lambda} } x^{-k-1} \theta_x \bullet H^k_n(x,\lambda),
\end{eqnarray}
where the second equality follows from $ \theta_{x} = x\frac{\pd{}}{\pd{}x} $ and~(\ref{eq:dxHkn_1}). 

\begin{remark}
\label{remark_change_variables}
Using the notation $  \varphi =\sqrt x$ and $  \psi=\sqrt \lambda $, with $ \psi $ constant, functions $ {}_0F_1(;n; x \lambda) $, $ H^k_n(x,\lambda) $, and $ v(x) $ become $ {}_0F_1(;n;  \varphi^2  \psi^2) $, $ H^k_n(  \varphi ^2, \psi^2) $, and $  v(  \varphi ) =e^{-2   \varphi   \psi} {}_0F_1(;n;  \varphi^2  \psi^2) $, respectively. Because $ \theta_\varphi = \varphi \frac{\pd{}}{\pd{}\varphi} = \frac{1}{2} \theta_x $ for $ x = \varphi^2 $, from~(\ref{equation_function_v}) we can write $ \theta_\varphi \bullet v(\varphi) = \frac{1}{2} \theta_x \left( e^x e^{-2\sqrt {x \lambda} } x^{-k-1} \right) \theta_x \bullet H^k_n(x,\lambda) $.
\end{remark} 

The following theorem gives the LODE system w.r.t.~$ \varphi $ satisfied by $ H^k_n(  \varphi ^2, \psi^2) $ and $ v(  \varphi ) $.

\begin{theorem}
\label{theorem: new-pfaffian}
The  vector valued function $\ug(  \varphi) = \begin{matrix} (H^k_n(  \varphi ^2, \psi^2) &  v(  \varphi ) & \theta_\varphi \bullet v(  \varphi ))^{\sf T}\end{matrix} $ satisfies the following $3$-dimensional LODE system:
\begin{eqnarray}
\label{equation_LODE_system_small}
\label{equation_LODE_system_small_R}
\frac{\pd{}}{\pd{}  \varphi } \ug (  \varphi) =\frac{1}{  \varphi } \left( 
\begin{array}{ccc}
0 & 2e^{-  \varphi ^2+2  \varphi   \psi}  \varphi ^{2(k+1)} & 0\\
0 & 0 & 1\\
0 & -2 (2n-1) \varphi   \psi & -\left[ 4   \varphi   \psi+2(n-1)\right] \\
\end{array}
\right) \ug (  \varphi).
\end{eqnarray}
\end{theorem}

\begin{IEEEproof}
See Appendix~\ref{section_proof_theorem_3d_LODE}; it uses $ {}_0F_1(;n;z) $ differential equation\cite[Eq.~(15.10.1)]{NIST_book_10}.
\end{IEEEproof}

\begin{theorem}  \label{th:stabile-LODE-2017-12-05}
The $3$-dimensional LODE system~(\ref{equation_LODE_system_small}) is stabile 
for $( \begin{matrix}H^k_n(\varphi^2,\psi^2) & v(\varphi) &\theta_\varphi \bullet v(\varphi) \end{matrix})^T$
when $\varphi \rightarrow \infty$, $ \forall \psi $.
\end{theorem}

\begin{IEEEproof}
Sketched in Appendix~\ref{proof_th:stabile-LODE-2017-12-05}.
\end{IEEEproof}

Next, we reveal that the stabile $3$-dimensional LODE system~(\ref{equation_LODE_system_small_R}) can be deduced systematically from the $4$-dimensional LODE system~(\ref{eq:pfaffian}), which is not stabile, and that this is possible in general.
Finally, Section~\ref{sec:numerical_results} demonstrates with numerical results that the HGM w.r.t.~$ \varphi $ based on~(\ref{equation_LODE_system_small}) yields reliable $ H^k_n(x,\lambda) $ computation even for very large values of $ \lambda $.

\subsection{Gauge Transformation-Based Deduction of Stabile LODEs}

Above, we derived the stabile LODE~(\ref{equation_LODE_system_small}) from the inhomogeneous differential equation for $H^k_n$ in~(\ref{eq:dxHkn_1}) after accounting for the asymptotic behavior of $ {}_0F_1(;\cdot;\cdot) $. Nevertheless, as we state below, in general, a stabile LODE can be found by applying a suitable gauge transformation to an available LODE. Gauge transformations are widely used in mathematical physics\cite{jimbo_np_81}. Our definition for them here is analogous but simpler.

\begin{definition}
Consider the LODE
\begin{equation} \label{eq:20171116-ode}
 {\bf f}'(x) = \frac{{\mathrm{d}}}{{\mathrm{d}} x}{\bf f}(x)= {\bf P}(x) {\bf f}(x),
\end{equation}
where ${\bf P}(x)$ is a matrix-valued function and 
${\bf f} (x)$ is a vector-valued function.
Let ${\bf G}(x)$ be an invertible matrix
and ${\bf h}(x)$ be the vector that satisfies 
${\bf f}={\bf G} {\bf h}$, so that 
$ {\bf f}'={\bf G}' {\bf h} + {\bf G} {\bf h}' $ and, hence,
\begin{equation} \label{eq:20171116-gauge-transformed}
 {\bf h}' = \left({\bf G}^{-1} {\bf P} {\bf G} -{\bf G}^{-1} {\bf G}' \right) {\bf h}.
\end{equation}
Then, ${\bf G}$ is called the gauge transformation that yields LODE~(\ref{eq:20171116-gauge-transformed}) from LODE~(\ref{eq:20171116-ode}).
\end{definition}
The proofs of the following two theorems are not included due to manuscript length limitation.

\begin{theorem} \label{th:by-gauge-2017-12-04}
A gauge transformation yields the stabile $ 3$-dimensional LODE system~(\ref{equation_LODE_system_small_R})
from the $ 4 $-dimensional LODE system~(\ref{eq:pfaffian}), which is not stabile, by block upper triangulation.
\end{theorem}

\begin{theorem}
\label{theorem_stabile_LODE_from_gauge}
From a given LODE system that is not stabile for a target function, a lower-dimensional stabile LODE system can be derived algorithmically by gauge transformations.
\end{theorem}

\begin{remark}
Gauge transformations can help convert LODEs obtained by hand and computer algebra into stabile LODEs that then guarantee the numerical convergence of HGM-based performance measure evaluation, also for other MIMO transceiver methods and fading models.
\end{remark}

\subsection{Further Enhancements to the Proposed HGM}
\label{section_enhanced_HGM}

Recall from Remark~\ref{remark_stabile_not} that if a LODE system is stabile for a function then it can be computed by HGM for arbitrary values of its variables, theoretically.
Practically,
we have found that HGM w.r.t.~$ x $ based on the stabile LODE system~(\ref{equation_LODE_system_small}) leads for large  $ \lambda $ to exceeding the largest number representable in \texttt{double} arithmetic in C\footnote{The largest \texttt{double} is $1.797693 \times 10^{308}$, whereas numerical integration in Mathematica yields $H^0_1(10^8,10^8) \approx 10^{43429447}$.}.
Nevertheless, we shall demonstrate that we can overcome this limitation by enhancing HGM with the following modifications:
\begin{enumerate}
\item Maintain values within bounds of the {\texttt{double}} representation by applying to the vector $ \ug $ from Theorem \ref{theorem: new-pfaffian} the following additional gauge transformations:
\begin{eqnarray}
\label{equation_gauge_G2}
\uG_2 & = &
\left( \begin{array}{ccc}
  \exp(-  \varphi  ^2+2  \varphi   \psi ) & 0 & 0 \\
   0 & 1 & 0 \\
   0 & 0 & 1 \\
\end{array} \right), \quad \text{for} \,   \varphi   <  \psi,
\end{eqnarray}
\begin{eqnarray}
\label{equation_gauge_G3}
\uG_3 & = & 
\left( \begin{array}{ccc}
  \exp( \psi ^2) & 0 & 0 \\
   0 & 1 & 0 \\
   0 & 0 & 1 \\
\end{array} \right), \quad \text{for} \,   \varphi   \geq  \psi.
\end{eqnarray}
On the one hand, in matrix $ \uG_2$, the exponent $-  \varphi  ^2+2  \varphi \psi $ takes its maximum value $ \psi ^2$ when $  \varphi  = \psi $ and becomes $ 0 $ when $  \varphi  = 0 $ or $  \varphi  =  2 \psi $.  On the other hand, matrix $ \uG_3$ is constant because $ \psi $ is considered constant. While these transformations have been selected heuristically, they have been found reliable.
For example, $x =  \varphi ^2 = 10^8$ and  $\lambda = \psi^2 = 10^8$ yield $ H^0_1(  \varphi ^2, \psi^2)/\exp(-  \varphi  ^2+2  \varphi   \psi) \approx 0.5 $, i.e., within the \texttt{double} representation range.  
\item 
Use numerical integration of~(\ref{equation_Hnk_definition_integral}) to obtain initial conditions, as series truncation no longer converges in a reasonable time for large  $ \lambda $.  
\item Use {\emph{big float}}\footnote{I.e., arbitrary precision. Note that, as required precision increases, computation time can also increase.} in computing the determinant of the matrix with elements $H^k_n(x, \lambda)$.  Relevant experiments have been implemented in Risa/Asir.
\item Estimate the numerical errors arising in the evaluation of the determinant by adding random errors of the same size as the numerical error of $H^k_n(x, \lambda)$ to each element of $\uPhi(x)$.
\end{enumerate}

Numerical results from HGM enhanced as above are shown in Sections~\ref{sec:numerical_large_lambda} and~\ref{section_MRC_results}.

\section{Numerical Results}
\label{sec:numerical_results}

Below, we describe our results from evaluating the CDF of the dominant eigenvalue $ \phi_s $ of $ (K + 1) \uHH \uH $ for several antenna number pairs $ (\NT, \NR) $ and $ \uSd $ eigenvalue sets $\ulambda=\{\lambda_1, \lambda_2, \cdots, \lambda_s\}$.
Recall that this CDF yields the MRC outage probability, after accounting for the scaling in~(\ref{equation_MRC_Outage_Probability}).
We shall find that, unlike conventional methods, our HGM is consistently reliable and efficient for CDF computation, even in extreme conditions, i.e., very large $ \lambda_i$'s (or $\NR$, $ \NT$, and $K$).

We employed a computer with the Debian ``wheezy'' operating system, a Xeon E5-4650 processor at 2.7 GHz, and 256 GB of memory.\label{computer_description}
Table~\ref{table_numerical_procedure} summarizes the numerical procedures and settings employed for CDF computation, whereas footnotes therein describe the software packages used for the implementation of these procedures as well as further settings. 

\begin{table}
\begin{center}
\caption{Numerical Procedures and Settings}
\label{table_numerical_procedure}
\begin{tabular}{l | l | l | l | l}
Value of $ \lambda $ & Section & HGM Initial Condition  &  HGM LODEs & Comparison \\ \hline \hline
Very small & \ref{experiment_HGM_wrt_lambda} & Truncation\footnote{\label{footnote_LODE_Integration}In C, with double precision representation.}~(\ref{eq:truncated}), $ \lambda_0 = 10^{-5} $  & Eq.~(\ref{eq:pfaffian_wrt_lambda}), w.r.t.~$\lambda$, \footnote{\label{footnote_LODE_only}In C, with Runge--Kutta method of step size $h=10^{-4}$.} &  Truncation\textsuperscript{\ref{footnote_LODE_Integration}}
~(\ref{eq:truncated}) \\ \hline
Small, moderate & \ref{sec:numerical_small_lambda} & Truncation\textsuperscript{\ref{footnote_LODE_Integration}}~(\ref{eq:truncated}), $ \lambda_0 = 10^{-5} $ & Eq.~(\ref{equation_LODE_system_small_R}), w.r.t.~$x$\textsuperscript{\ref{footnote_LODE_only}} & Integration\footnote{In Mathematica, to 16 significant digits (by default).} , Eq.~(\ref{eq:integral-rep}) \\ \hline
Very large & \ref{sec:numerical_large_lambda} & Integration\footnote{\label{footnote_Math_Integration}In Mathematica, to 20 significant digits.}, Eq.~(\ref{eq:integral-rep}) & Eq.~(\ref{equation_LODE_system_small_R}), w.r.t.~$x$\footnote{\label{footnote_Celeron} HGM enhanced as described in Section~\ref{section_enhanced_HGM}, in C, with adaptive Runge--Kutta method, with absolute error $ 10^{-5}$ and relative error $10^{-20}$.
The significand and exponent of big float are set to $7$ digits and $32$ bits, respectively. 
Values of $H^k_n(x,\lambda)$ stored every time $ x $ increased by $10^4$ in the interval $ [1 \times 10^6, 2.6 \times 10^8]$.} & Simulation\footnote{\label{footnote_simulation_R}In the statistical computing system R along with the package \texttt{cmvnorm} comprising various utilities for the complex multivariate Gaussian distribution\cite{software_cmvnorm_17}, Monte-Carlo simulation of $ 10^6 $ random $ 5 \times 7 $ matrix samples.}  for $x=2 \times 10^8$ \\ \hline
Large, very large & \ref{section_MRC_results} & Integration\textsuperscript{\ref{footnote_Math_Integration}}, Eq.~(\ref{eq:integral-rep}) & Eq.~(\ref{equation_LODE_system_small_R}), w.r.t.~$x$\textsuperscript{\ref{footnote_Celeron}}& Simulation\footnote{In Matlab.}, Fig.~\ref{figure_Outage_Probability_vs_Gb_HGM_Sim_NT_5_NR_5_Takayama} only. \\ \hline
\end{tabular}
\end{center}
\end{table}

\subsection{Numerical Results for Small-to-Moderate $ \NT$, $ \NR $, and $ \lambda_i $, from HGM w.r.t.~$ x $ and Integration}
\label{sec:numerical_small_lambda}

\subsubsection{Numerical Evaluation Procedure}
Here, we used unenhanced HGM w.r.t.~$x$ based on the LODE system~(\ref{equation_LODE_system_small}) to evaluate $H^k_n(x,\lambda)$. 
For HGM, we computed the initial condition
$\ug( \varphi _0,\psi)$
using the series truncation~(\ref{eq:truncated}) for a  value of $ \varphi _0$ sufficiently small to ensure numerical convergence, and then we computed $\ug( \varphi ,\psi)$ based on~(\ref{equation_LODE_system_small}) using the Runge--Kutta method.
We compared with Mathematica numerical integration as in~(\ref{eq:integral-rep}), because series truncation~(\ref{eq:truncated}) does not always converge numerically. 
Table~\ref{table_numerical_procedure}, third line, details these procedures and settings.
Finally, substituting $H^k_n(x,\lambda)$ in CDF expression~(\ref{equation_CDF_det_original}), yielded $ \pr(\phi_s \leq x) $ for $ \NT = 5 $, $ \NR = 5 : 9 $, and for two sets of eigenvalues $ \ulambda = \{ \lambda_1, \cdots, \lambda_s \} $ for $ \uSd  $, as described below.

\subsubsection{Results for $ \NT = 5 $, $ \NR = 5 : 9 $, and Small $ \lambda_i $}
\label{section_few_antennas_small_K}
Here, we compare the accuracy and computation time of HGM and numerical integration. We let $(\NT,\NR)=(5,5:9)$, to validate the HGM for small antenna numbers and to observe the effect of $ \NR $, and let the set of eigenvalues for $ \uSd $ be $\ulambda=\{0.1, 0.2, 0.3, 0.4, 0.5\}$, i.e., small.

\begin{figure}
\begin{center}
\includegraphics[width=3.45in]
{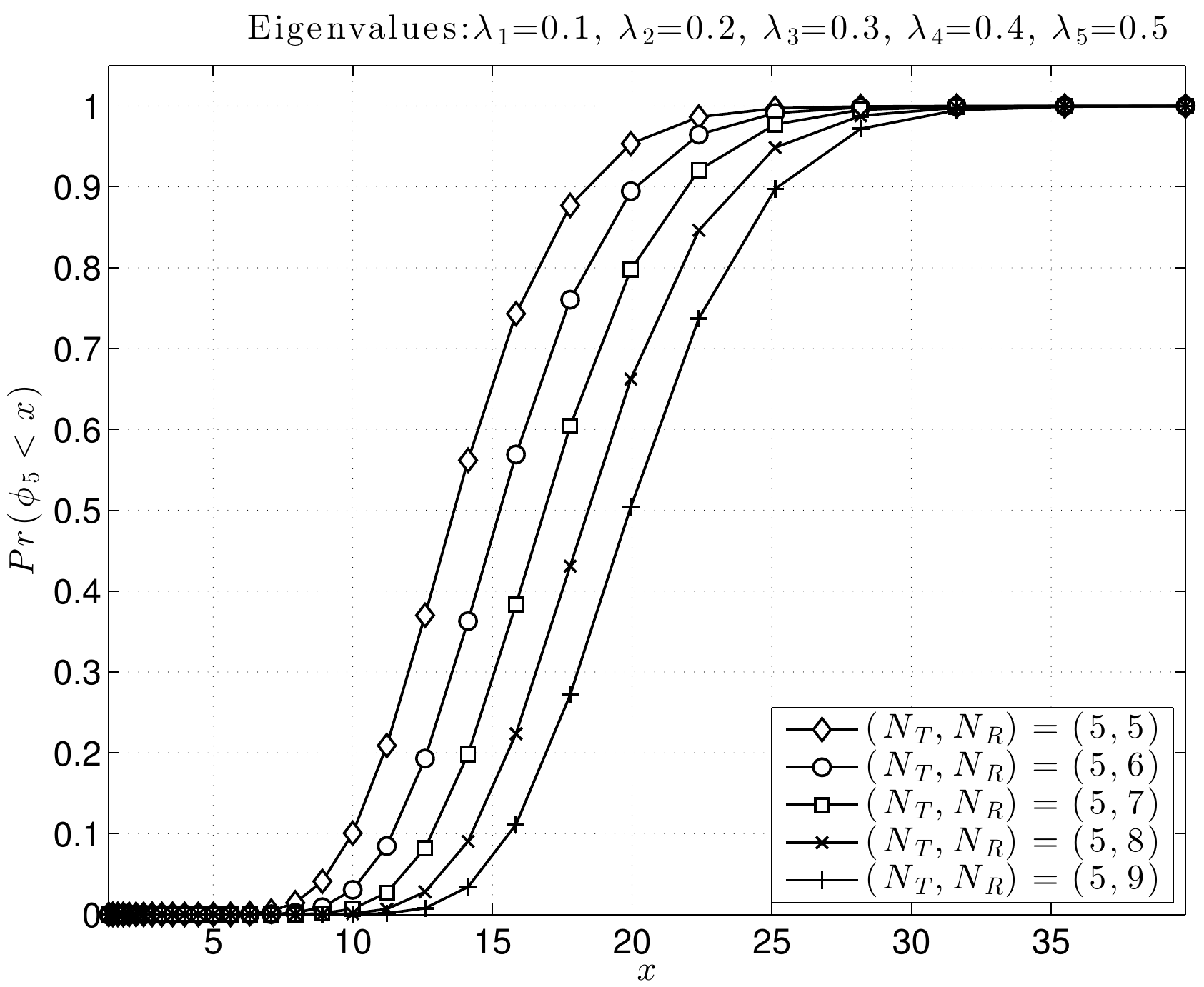}
\caption{CDF of dominant eigenvalue, $ \pr(\phi_5 \leq x) $, from HGM, for $ \NT = 5 $, $ \NR = 5 : 9 $, and the set of eigenvalues $\ulambda=\{0.1, 0.2, 0.3, 0.4, 0.5\}$ for $ \uSd = (K + 1) \uHdH \uHd $.}
\label{figure_Eigenvalue_CDF_vs_x_NT_5_NR_5_9_Fadil}
\end{center}
\end{figure}

Fig.~\ref{figure_Eigenvalue_CDF_vs_x_NT_5_NR_5_9_Fadil} shows the CDF of the dominant eigenvalue, i.e., $ \pr(\phi_5 \leq x) $, computed with HGM for the mentioned $ (\NT, \NR) $ pairs, confirming that the CDF grows from $ 0 $ to $ 1 $. It also reveals that, given $ x $, a larger $ \NR $ yields lower $ \pr(\phi_5 \leq x) $, i.e., lower outage probability for MRC, confirming that employing more antennas yields better performance.

Table~\ref{Table_results_NT_5_NR_5_8} reveals that HGM is computationally much less demanding than Mathematica numerical integration.  Increasing $ \NR $ slightly increases the computation time for both methods. 

\renewcommand{\arraystretch}{0.8}

\begin{table}
    \begin{center}
\caption{Computation time, in seconds, for Section~\ref{section_few_antennas_small_K}}
\label{Table_results_NT_5_NR_5_8}
        \begin{tabular}{c|l|l}
            \hline
            {$(\NT,\NR)$}  & {{Integration}} & {HGM} \\
            \hline
            (5,5)    &    3.28     & 0.244 \\
            (5,6)    &    3.59     & 0.256 \\
            (5,7)    &    3.713    & 0.260 \\
            (5,8)    &    3.834    & 0.244 \\
            (5,9)    &    3.87     & 0.256 \\
            \hline
    \end{tabular}
    \end{center}
\end{table}

\subsubsection{Results for Moderate $ \NT $, $ \NR $, and $ \lambda_i $}
\label{section_moderate_antennas_K}

Compared to above, herein, we consider a case with more antennas and larger eigenvalues, i.e., $ (\NT, \NR ) = (10, 10) $ and $\ulambda=\{1,2,\ldots,10\}$ as eigenvalue set for $ \uSd = (K + 1) \uHdH \uHd $. 
Based on Remark~\ref{remark_trace_KNTNR_1}, these choices yield $ K \approx - 2.2 $~dB, i.e., a small value.
Table~\ref{table_evaluation_lambda_set_1} shows good agreement between the CDF evaluations by HGM and Mathematica integration.
However, we have found that HGM required only $ 0.920 $ seconds vs.~$7.99 $ seconds required by Mathematica integration.

\begin{table}
\begin{center}
    \caption{$\pr(\phi_{10} \leq{}x)$ for  $ (\NT, \NR ) = (10, 10) $, $ \ulambda_1=\{1,2,\ldots,10\}$, for Section~\ref{section_moderate_antennas_K}}
    \label{table_evaluation_lambda_set_1}
    \begin{tabular}{c | l | l}
        \hline
        $ \log_{10} x $ & Integration & HGM\\
        \hline
        1 & 9.52818e-45    & 3.77131e-30\\
        1.1 & 3.6096e-30    & 5.64982e-23\\
        1.2 & 1.58216e-16    & 1.56884e-16\\
        1.3 & 5.18699e-11        & 5.21756e-11\\
        1.4 & 1.36179e-06    & 1.36963e-06\\
        1.5    & 0.00203227    & 0.0020352\\
        1.6    & 0.148478    & 0.14857\\
        1.7    & 0.781498    & 0.781594\\
        1.8    & 0.99524    & 0.995231\\
        1.9 & 1.00001    & 0.999994    \\
        2        & 1.00001    & 0.999998\\
        \hline
    \end{tabular}    
\end{center}
\end{table}

These results demonstrate that the proposed HGM w.r.t.~$x$ based on the stabile LODE system~(\ref{equation_LODE_system_small}) is computationally reliable and efficient for evaluating the performance of moderately-sized MIMO systems. 
However, because the very small value ensuing for $ K $ does not characterize realistic Rician fading, we next attempt the CDF calculation for larger eigenvalues for $ \uSd = (K + 1) \uHdH \uHd $, and, thus, larger values for $ K $.
We have found that such larger values have required the application of HGM w.r.t.~$x$ enhanced as described in Section~\ref{section_enhanced_HGM}.

\subsection{Numerical Results for Very Large $ \lambda_i $, from Enhanced HGM w.r.t.~$ x $}
\label{sec:numerical_large_lambda}

These numerical experiments have been implemented as described in Table~\ref{table_numerical_procedure}, fourth line, for $ \NT = 5 $, $ \NR = 5:7$, and $ \ulambda=\{0.4, 0.8, 1.2, 1.6, 2 \}\times10^8 $.
Table~\ref{table:large_eigen_outage_result_3} shows the CDF output by our enhanced HGM, confirming that it increases from 0 to $1$ with increasing $x$. However, Table~\ref{table:large_eigen_comp_time_2} indicates that HGM computation time is relatively large.
On the other hand, series truncation~(\ref{eq:truncated}) and Mathematica integration both fail for such large eigenvalues for $ \uSd $.
\begin{table}
\begin{center}
\caption{CDF $\pr(\phi_5 \leq{}x)$ output by HGM, for Section~\ref{sec:numerical_large_lambda}}
\label{table:large_eigen_outage_result_3}
\begin{tabular}{c | l | l | l}
 $x$        & $\pr(\phi_5 \leq{}x)$, $5\times 5$  &  $\pr(\phi_5 \leq{}x)$, $5\times 6$ & $\pr(\phi_5 \leq{}x)$, $5\times 7$\\ \hline
$1.9990 \times 10^8$ &  2.841503e-07 &  2.840696e-07 &  2.840035e-07 \\ 
$1.9991 \times 10^8$ &  3.372661e-06 &  3.371785e-06 &  3.371082e-06 \\ 
$1.9992 \times 10^8$ &  3.147594e-05 &  3.146851e-05 &  3.146270e-05 \\ 
$1.9993 \times 10^8$ &  0.00023143703 &  0.00023138784 &  0.00023135066 \\ 
$1.9994 \times 10^8$ &  0.0013442040 &  0.0013439501 &  0.0013437654 \\ 
$1.9995 \times 10^8$ &  0.0061883532 &  0.0061873284 &  0.0061866184 \\ 
$1.9996 \times 10^8$ &  0.022687564 &  0.022684307 &  0.022682226 \\ 
$1.9997 \times 10^8$ &  0.066662879 &  0.066654839 &  0.066650189 \\ 
$1.9998 \times 10^8$ &  0.15839370 &  0.15837759 &  0.15836978 \\ 
$1.9999 \times 10^8$ &  0.30816435 &  0.30813940 &  0.30812954 \\ \hline
$2.0000 \times 10^8$ &  0.49958230 &  0.49954954 &  0.49954438 \\ \hline
$2.0001 \times 10^8$ &  0.69109536 &  0.69106073 &  0.69105953 \\ 
$2.0002 \times 10^8$ &  0.84109309 &  0.84106275 &  0.84107198 \\ 
$2.0003 \times 10^8$ &  0.93306099 &  0.93303238 &  0.93305115 \\ 
$3.0000 \times 10^8$ &  1.000017 &  0.99999227 &  1.000260
\end{tabular}
\end{center}
\end{table}

\renewcommand{\arraystretch}{0.9}

\begin{table}
\begin{center}
\caption{HGM computation time, in seconds and hours, for Section~\ref{sec:numerical_large_lambda}}
\label{table:large_eigen_comp_time_2}
\begin{tabular}{c|c|c|c}
$\NT \times \NR$& Initial condition & LODE solving & HGM Total \\ \hline
$5\times 5$ & $ 30 $ s  & 8,037 s & 8,067 s $ \approx $ 2.25 h \\
$5\times 6$ & $ 41 $ s  & 25,340 s & 25,381 s $ \approx $ 7 h \\ \hline
$5\times 7$ & $ 45 $ s & 50,099 s & 50,144 s $ \approx $ 14 h \\ \hline
\end{tabular}
\end{center}
\end{table}

For HGM validation, we did a Monte-Carlo simulation of $ 10^6 $ random $ 5 \times 7 $ matrix samples, as shown in Footnote~\ref{footnote_simulation_R} on page~\pageref{footnote_simulation_R}. This simulation output the value $0.499458$ for $\pr(\phi_5 \leq{} x=2 \times 10^8)$ whereas HGM output $ 0.49954438 $ (see Table~\ref{table:large_eigen_outage_result_3}). The simulation took $ 290 $~s whereas HGM required $ 50144 $~s (see Table~\ref{table:large_eigen_comp_time_2}).
Nevertheless, the Runge--Kutta procedure in HGM inherently computes the CDF at numerous samples throughout the range of interest for $ x $.

\subsection{MIMO MRC Outage Probability Results for Large $ K $ and $ \NR $, from Enhanced HGM w.r.t.~$ x $}
\label{section_MRC_results}
These numerical experiments have been implemented as described in Table~\ref{table_numerical_procedure}, last line. 

\subsubsection{MRC Outage Probability Results for $ \NT = 5 $, $ \NR = 5$, and Very Large $ K $}

\begin{figure}
\begin{center}
\includegraphics[width=3.45in]
{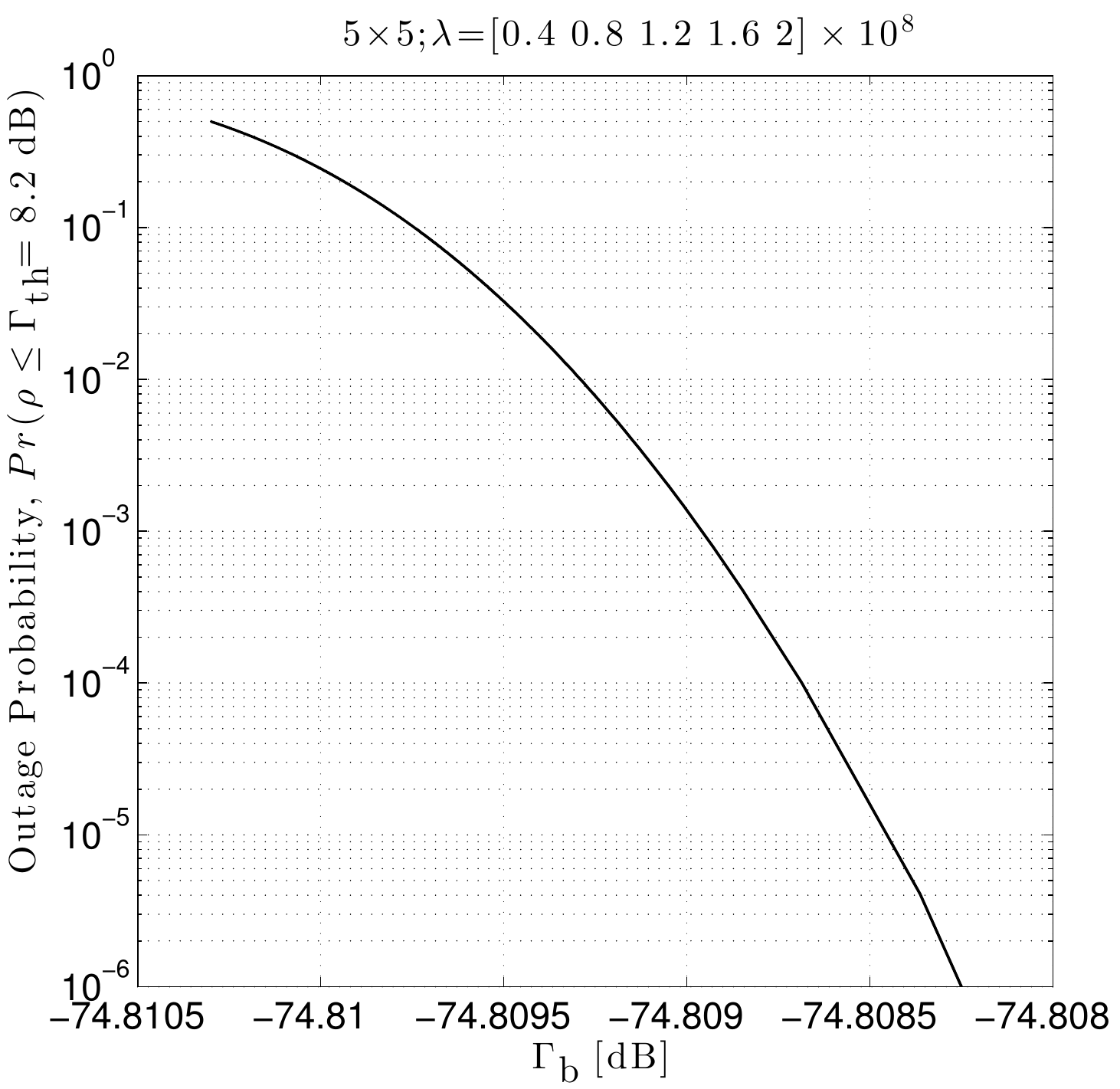}
\caption{Outage probability $ \pr(\rho_{\text{MRC}} \leq \Gamma_{\text{th}} = 8.2 \text{ dB}) = \pr \left(\phi_s \leq \frac{(K+1) \Gamma_{\text{th}}}{\Gamma_{\text{b}}} \right)  $ vs.~$ \Gamma_{\text{b}} $, from HGM, for $ \NT = 5 $, $ \NR = 5$, and set of eigenvalues $\ulambda=\{0.4, 0.8, 1.2, 1.6, 2 \} \times 10^8$.}
\label{figure_Outage_Probability_vs_Gb_NT_5_NR_5_Takayama}
\end{center}
\end{figure}

\begin{figure}
\begin{center}
\includegraphics[width=3.45in]
{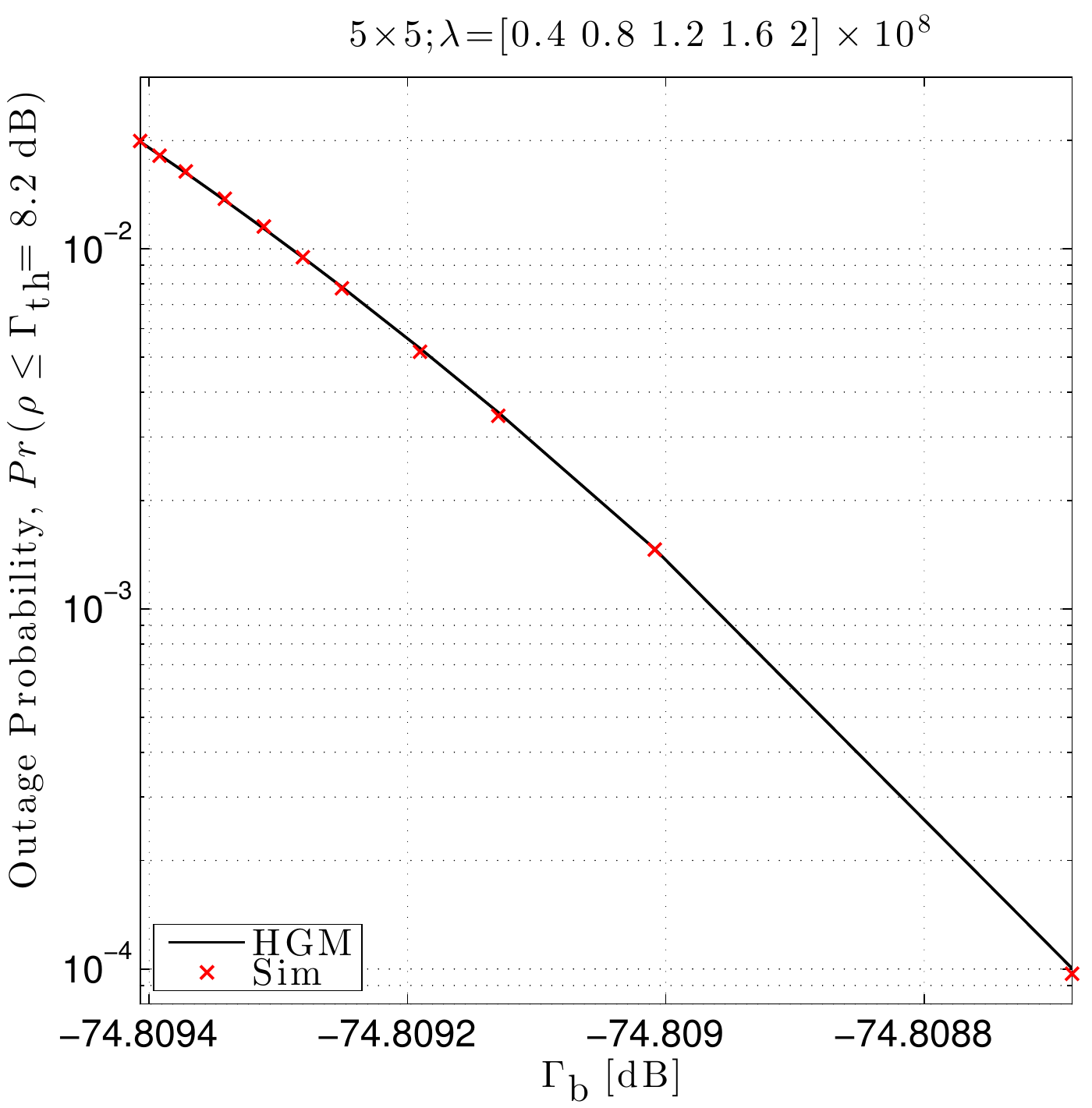}
\caption{Outage probability $ \pr(\rho_{\text{MRC}} \leq \Gamma_{\text{th}} = 8.2 \text{ dB}) = \pr \left(\phi_s \leq \frac{(K+1) \Gamma_{\text{th}}}{\Gamma_{\text{b}}} \right)  $ vs.~$ \Gamma_{\text{b}} $, from HGM and simulation, for $ \NT = 5 $, $ \NR = 5$, and set of eigenvalues $\ulambda=\{0.4, 0.8, 1.2, 1.6, 2 \} \times 10^8$.}
\label{figure_Outage_Probability_vs_Gb_HGM_Sim_NT_5_NR_5_Takayama}
\end{center}
\end{figure}

Fig.~\ref{figure_Outage_Probability_vs_Gb_NT_5_NR_5_Takayama} shows the outage probability $ \pr(\rho_{\text{MRC}} \leq  8.2 \text{ dB}) $ vs.~$ \Gamma_{\text{b}} $ from HGM for $ \NT = 5 $, $ \NR = 5 $, and eigenvalue set $ \ulambda=\{ (2/5) , (4/5), (6/5), (8/5), 2 \} \times  10^8$ for $ \uSd $. Then, from Remark~\ref{remark_trace_KNTNR_1}, we have that $ \sum_{i = 1}^{5} \lambda_i = 6 \times 10^8 $ yields $ K \approx 70 $~dB.
Fig.~\ref{figure_Outage_Probability_vs_Gb_HGM_Sim_NT_5_NR_5_Takayama} validates the HGM results against Matlab simulations, demonstrating that our enhanced HGM is reliable for the evaluation of MRC under Rician fading.
Additionally, the HGM efficiency has allowed for a wider outage probability range to be explored vs.~simulations --- compare the vertical ranges in Figs.~\ref{figure_Outage_Probability_vs_Gb_NT_5_NR_5_Takayama} and~\ref{figure_Outage_Probability_vs_Gb_HGM_Sim_NT_5_NR_5_Takayama}.

\subsubsection{MRC Outage Probability Results for $ \NT = 5 $, $ \NR = 100 $, and $ K \approx 20 $~dB}

\begin{figure}
\begin{center}
\includegraphics[width=3.45in]
{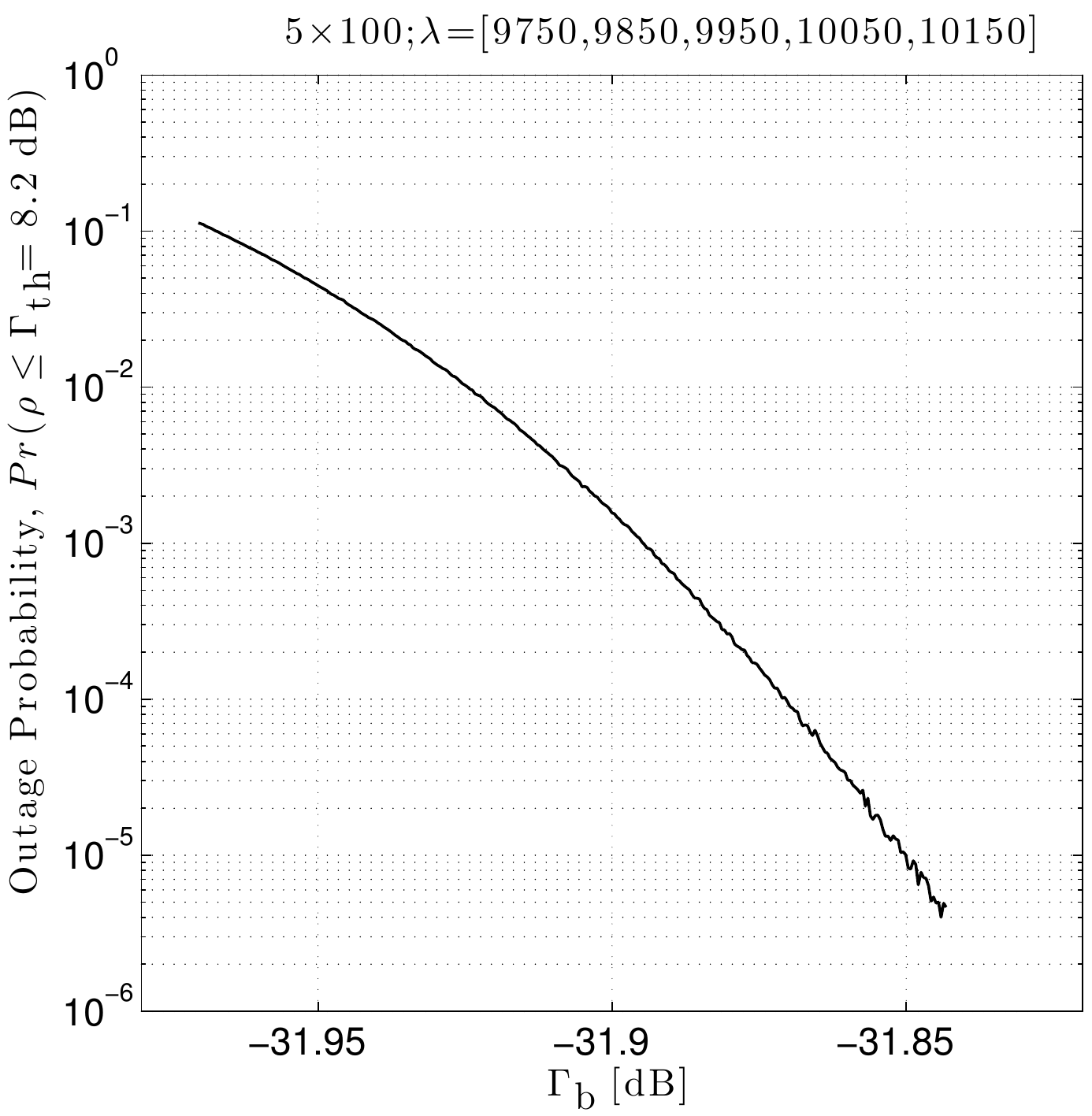}
\caption{Outage probability $ \pr(\rho_{\text{MRC}} \leq \Gamma_{\text{th}} = 8.2 \text{dB}) = \pr \left(\phi_s \leq \frac{(K+1) \Gamma_{\text{th}}}{\Gamma_{\text{b}}} \right)  $ vs.~$ \Gamma_{\text{b}} $, obtained by HGM for $ \NT = 5 $, $ \NR = 100 $, and the set of eigenvalues $ \ulambda=\{ 9750, 9850, 9950, 10050, 10150 \} $, i.e., $ \uSd $ condition number $ \lambda_5/\lambda_1=1.0410 $.}
\label{figure_Outage_Probability_vs_Gb_NT_5_NR_100_Takayama}
\end{center}
\end{figure}

\begin{figure}
\begin{center}
\includegraphics[width=3.45in]
{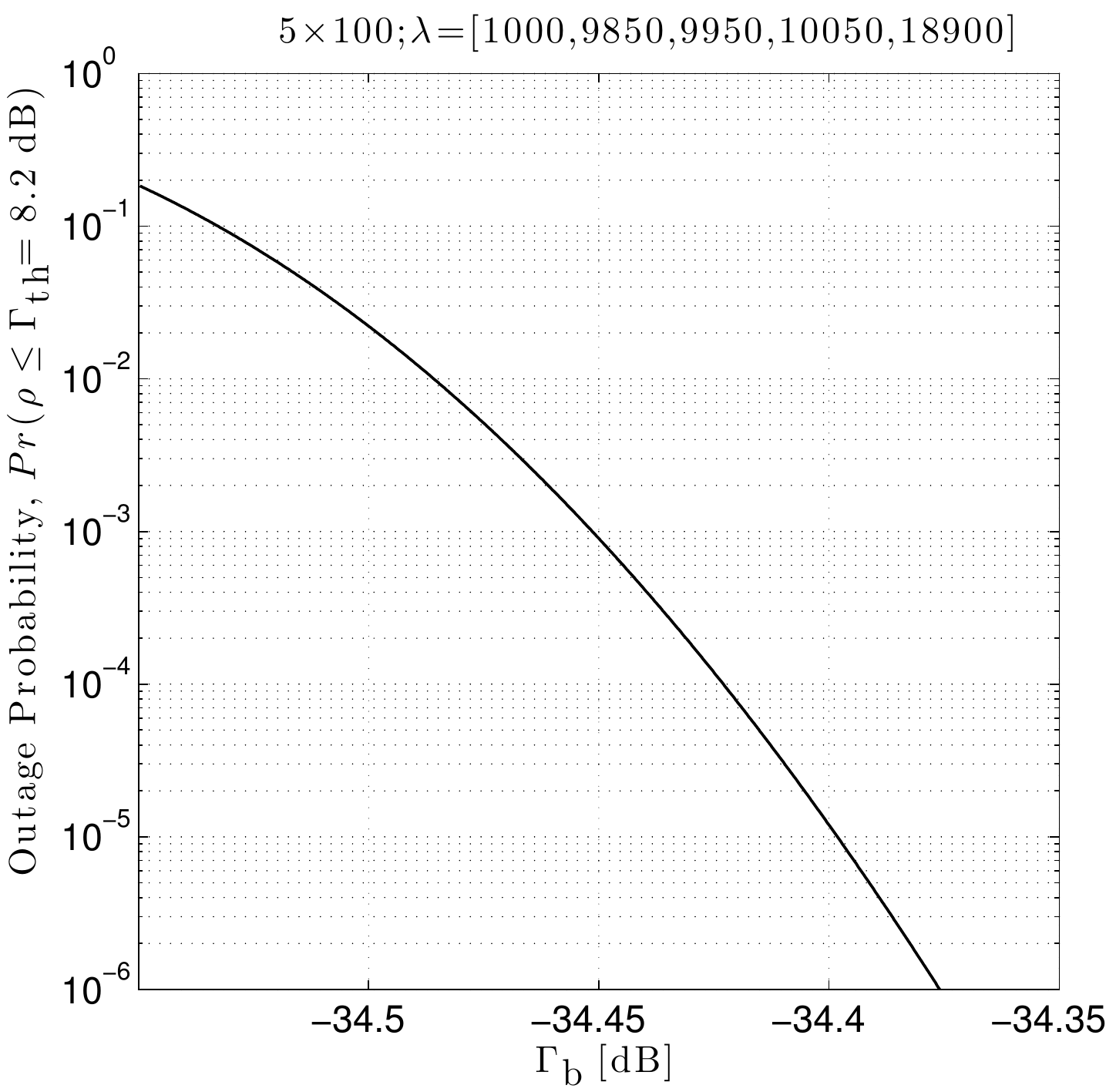}
\caption{Outage probability $ \pr(\rho_{\text{MRC}} \leq \Gamma_{\text{th}} = 8.2 \text{ dB}) = \pr \left(\phi_s \leq \frac{(K+1) \Gamma_{\text{th}}}{\Gamma_{\text{b}}} \right)  $ vs.~$ \Gamma_{\text{b}} $, obtained by HGM for $ \NT = 5 $, $ \NR = 100 $, and the set of eigenvalues $ \ulambda=\{ 1000,9850,9950,10050,18900 \} $, i.e., $ \uSd $ condition number $ \lambda_5/\lambda_1 = 18.9 $.}
\label{figure_Outage_Probability_vs_Gb_NT_5_NR_100_Takayama_Large_CN}
\end{center}
\end{figure}

Figs.~\ref{figure_Outage_Probability_vs_Gb_NT_5_NR_100_Takayama} and \ref{figure_Outage_Probability_vs_Gb_NT_5_NR_100_Takayama_Large_CN} show the outage probability $ \pr(\rho_{\text{MRC}} \leq8.2 \text{ dB}) $  vs.~$ \Gamma_{\text{b}} $ from our enhanced HGM, for $ \NT = 5 $, $ \NR = 100$, for $ \ulambda=\{ 9750, 9850, 9950, 10050, 10150 \} $ and $ \ulambda=\{ 1000,9850,9950,10050,18900 \} $, respectively.
Note that these eigenvalue sets yield the same $ \tr(\uSd) = \sum_{i=1}^{s} \lambda_i = \NT \NR K = 49750 $, i.e., the results in Figs.~\ref{figure_Outage_Probability_vs_Gb_NT_5_NR_100_Takayama} and \ref{figure_Outage_Probability_vs_Gb_NT_5_NR_100_Takayama_Large_CN} are for $ K \approx 20 $~dB.
On the other hand, these eigenvalue choices yield for the condition number $ \lambda_5/\lambda_1$ values of $ 1.0410 $ and $ 18.9 $, respectively.
Comparing Figs.~\ref{figure_Outage_Probability_vs_Gb_NT_5_NR_100_Takayama} and \ref{figure_Outage_Probability_vs_Gb_NT_5_NR_100_Takayama_Large_CN} reveals that a larger condition number for $ \uSd $ yields better MRC performance\footnote{This is because: 1) large $ K $ yields $ \uH \approx \uHd $, i.e., the largest eigenvalue $ \phi_s $ of $ (K+1) \uHH \uH $ approaches the largest eigenvalue of $ (K+1)\uHdH \uHd = \uSd $, and 2) the largest eigenvalue of $ \uSd $ increases, given $ \text{tr}(\uSd) $, when $ \text{rank}(\uSd) $ approaches $ 1 $, i.e., for increasing condition number for $ \uSd $.}.

\section{Summary and Conclusions}

This paper has sought the reliable CDF evaluation for the largest eigenvalue of complex noncentral Wishart matrices. A known CDF expression involves a determinant whose elements are integrals in a hypergeometric function. For these elements we have derived an infinite series and the satisfied LODEs, by hand and by using computer algebra to reveal underlying Gr\"obner bases. Our ensuing HGM deployment has been found to suffer from numerical instability. 
Therefore, we have also derived by hand a LODE system which we have proved stabile for its solution, i.e., theoretically ensuring numerical stability. 
Importantly,  we have revealed that stabile LODEs can be obtained algorithmically by gauge transformations of other LODEs, {\underline{in general}}.
Thus, unlike numerical truncation and integration, the ensuing HGM, enhanced with several gauge transformations and other modifications, has been proved reliable for parameter values as large as may occur in practice for massive MIMO system under Rician fading. 
Our enhanced HGM has also been demonstrated to be computationally more efficient than numerical integration.

\appendices

\section{Proof of Theorem~\ref{theorem_4d_LODE_pfaffians}}
\label{appendix_proof_theorem_ODE4}

\subsection{Obtaining Differential Equations for $H^k_n(x,\lambda)$ Based on Infinite Series~(\ref{eq:series-rep}), by Hand}

The proof follows from the fact, proven below, that $H^k_n(x,\lambda)$ satisfies the differential equations
\begin{eqnarray}
\label{eq:diff-eq}
\{\theta_{\lambda}(\theta_{\lambda}+n-1)+\lambda (\theta_{x}-\theta_{\lambda}-k-1)\}\bullet  H^k_n(x,\lambda) & = & 0, \\
\label{eq:diff-eq_1}
\{\theta_{x}(\theta_{x}-\theta_{\lambda}-k-1)+x\theta_{x}\}\bullet  H^k_n(x,\lambda) & = & 0.
\end{eqnarray}
First, according to the new infinite series for $ H^k_n(x,\lambda)  $ from~(\ref{eq:series-rep}), let us define the function
$f(x_1,x_2) = \sum_{p=0}^\infty  \sum_{q=0}^\infty c_{pq} x_1^p x_2^q$, 
where $  c_{pq} = \frac{(-1)^q (k+1)_{p+q}}{(n)_p (1)_p (1)_q (k+2)_{p+q}}$, $x_1=x\lambda$, $x_2 =x$.
It can be shown that $c_{pq}$ satisfies the following recurrence relations:
\begin{eqnarray}
\label{equation_c_recurrence_1}
(n+p)(1+p)(k+2+p+q) c_{p+1,q} &=& (k+1+p+q) c_{pq}, \\
\label{equation_c_recurrence_2}
(1+q)(k+2+p+q) c_{p,q+1} &=& -(k+1+p+q) c_{pq}.
\end{eqnarray}
Then, let us define the differential operators $ \theta_i = x_i \frac{\partial}{\partial x_i}$, $ i = 1, 2 $. Note that $\theta_i^m \bullet x_i^p = p^m x_i^p$. Thus, for any polynomial $\ell(\theta_1,\theta_2)$
we have
\begin{eqnarray}
\label{equation_l_x1_x2}
\ell(\theta_1,\theta_2) \bullet x_1^p x_2^q = \ell(p,q) x_1^p x_2^q.
\end{eqnarray}
From~(\ref{equation_c_recurrence_1}) and~(\ref{equation_l_x1_x2}), we have
\begin{eqnarray} (n+\theta_1-1) \theta_1 (k+1+\theta_1+\theta_2) \bullet c_{p+1,q} x_1^{p+1} x_2^q
= x_1 (k+1+\theta_1+\theta_2) \bullet c_{pq} x_1^p x_2^q,
\end{eqnarray}
which implies
\begin{eqnarray}
\label{equation_operator_theta1_theta2_1}
[(n+\theta_1-1) \theta_1 (k+1+\theta_1+\theta_2) - x_1 (k+1+\theta_1+\theta_2)] \bullet
f(x_1,x_2) = 0.
\end{eqnarray}
Analogously, from~(\ref{equation_c_recurrence_2}) and~(\ref{equation_l_x1_x2}), we can show that
\begin{eqnarray}
\label{equation_operator_theta1_theta2_2}
[\theta_2 (k+1+\theta_1+\theta_2) + x_2(k+1+\theta_1+\theta_2)] \bullet f(x_1,x_2) = 0.
\end{eqnarray}
Differential equations~(\ref{equation_operator_theta1_theta2_1}) and~(\ref{equation_operator_theta1_theta2_2}) can be recast in terms of  $x=x_2$ and $\lambda=x_1/x_2$, after using the chain rule of differentiation to obtain
$\theta_1 = \theta_\lambda$ and $\theta_2 = \theta_x-\theta_\lambda$, as follows:
\begin{eqnarray}
&& [(n+\theta_\lambda-1) \theta_\lambda (k+1+\theta_x) - x \lambda (k+1+\theta_x)] \bullet
f(x\lambda,x) = 0 \label{eq:eee1}, \\
&& [(\theta_x-\theta_\lambda) (k+1+\theta_x) + x (k+1+\theta_x)] \bullet f = 0.
\label{eq:eee2}
\end{eqnarray}
Multiplying (\ref{eq:eee2}) by $\lambda$ from the left and adding it to 
(\ref{eq:eee1}), we obtain:
\begin{eqnarray}
\label{equation_final_1}
(k+1+\theta_x) 
[(n+\theta_\lambda-1) \theta_\lambda + \lambda (\theta_x-\theta_\lambda) ]
\bullet f(x\lambda,x) = 0.
\end{eqnarray}
On the other hand, it can be shown by direct calculation that the following holds:
\begin{eqnarray}
\label{equation_final_2}
[(n+\theta_\lambda-1) \theta_\lambda + \lambda (\theta_x-\theta_\lambda) ] \bullet
f(x\lambda,x) = 0.
\end{eqnarray}
Finally, from\cite[Prop.~6.1.3]{hibi_book_13}  we have identity
$\ell(\theta_x,\theta_\lambda) x^k = x^k \ell(\theta_x+k,\theta_\lambda)$, which implies that
$x^k \ell(\theta_x,\theta_\lambda) \bullet \hat{f}(x, \lambda)
=\ell(\theta_x-k,\theta_\lambda) \bullet x^k \hat{f}(x, \lambda)
$ 
for any smooth function $ \hat{f}(x, \lambda) $. Using this property along with~(\ref{equation_final_1}) and~(\ref{equation_final_2}) yields~(\ref{eq:diff-eq}) and~(\ref{eq:diff-eq_1}), respectively.

Inconveniently, the differential equations in~(\ref{equation_final_1}) and~(\ref{equation_final_2}) mix the differential operators $ \theta_x $ and $ \theta_\lambda$. 
Therefore, we next attempt to obtain LODEs for $H^k_n(x,\lambda)$ w.r.t.~$x $ and $ \lambda $.  

\subsection{Obtaining LODEs for $H^k_n(x,\lambda)$ Using Gr\"obner Basis Computation, by Computer Algebra}
\label{section_Grobner_bases}

In general, given some  differential operators that annihilate a function, we can obtain LODEs for this function by computing the Gr\"obner basis of a left ideal on the ring of the available differential operators\cite[Ch.~6]{hibi_book_13}. 

Thus, let $ {\mathcal{R}} =\C (x,\lambda) \langle \frac{\pd{}}{\pd{}x}, \frac{\pd{}}{ \pd{}\lambda }\rangle$ be the ring of differential operators whose coefficients are complex-valued rational functions of $ (x,\lambda) $, and define the following bivariate operators which arise from~(\ref{eq:diff-eq}) and~(\ref{eq:diff-eq_1}), respectively,
\begin{eqnarray}
\label{equation_l1}
l_1 &=& \theta_{\lambda}(\theta_{\lambda}+n-1)+\lambda (\theta_{x}-\theta_{\lambda}-k-1),\\
\label{equation_l2}
l_2 &=& \theta_{x}(\theta_{x}-\theta_{\lambda}-k-1)+x\theta_{x},
\end{eqnarray} 
and which are elements of ring $ {\mathcal{R}} $.  Denote by $ {\mathcal{I}} $ the left ideal of $ {\mathcal{R}} $ generated by $l_1$ and $l_2$.  
We can now compute Gr\"obner bases for $ {\mathcal{I}} $, e.g., by using computer algebra systems; we have used Risa/Asir\cite{Risa_software_17}.

For our two-variable ring $ {\mathcal{R}} $, the computed Gr\"obner basis depends on the selected monomial ordering\cite[pp.~2, 7]{hibi_book_13}\cite[Sec.~2.4]{buchberger_mssp_01}:
\begin{itemize}
\item
Lexicographic ordering with $\frac{\partial}{\partial x} < \frac{\partial}{\partial \lambda}$, which yields the Gr\"obner basis $ {\mathcal{G}} = \{l_3,l_4,l_5\}$: 
\begin{eqnarray}\label{basis-x}
\label{equation_l3}
l_3 &=& \left[ (\theta_x+x-k-1)(\theta_x+x-k+n-2)-x \lambda \right] \theta_x,\\
l_4 &=& (\theta_{\lambda} - \theta_x - x+k+1)\theta_x, \nonumber \\
l_5 &=& \theta_{\lambda}^2 -(\lambda-n+1)\theta_{\lambda} + \lambda \theta_x -  \lambda(k+1), \nonumber
\end{eqnarray}
Note also that this basis implies that the holonomic rank of the system is $4$.
\item
Lexicographic ordering with $\frac{\partial}{\partial \lambda} < \frac{\partial}{\partial x}$, which yields the Gr\"obner basis $ {\mathcal{G}}' = \{l_6,l_7\}$:
\begin{eqnarray}
\label{basis-lambda}
l_6 & = & \theta_{\lambda}^4+(\lambda+2n-4)\theta_{\lambda}^3+\left[ n^2-5n+5-(x+k+n) \lambda \right]\,\theta_{\lambda}^2  +\left[ x\lambda^2-(n-1)x \lambda \right. \nonumber \\
&& \left. -(k+1)(n-1) \lambda -(n-2)(n-1)\right] \theta_{\lambda} +(k \lambda+1)x \lambda, \label{equation_l6} \\
l_7 & = &\theta_x+\frac{1}{\lambda}\theta_{\lambda}^2+\left(-1+\frac{n-1}{\lambda}\right)\theta_{\lambda}-(k+1). \nonumber
\end{eqnarray}
Note that $\theta_x$ does not appear in the differential operator $l_6$.
\end{itemize}

Thus, by using lexicographic ordering, our computer-algebra-based derivation of the Gr\"obner bases for the bivariate differential operators $ l_1 $ and $ l_2 $ yielded the univariate operators from~(\ref{equation_l3}) and~(\ref{equation_l6}), which annihilate function $ H^k_n(x,\lambda) $. They can be recast as the Theorem~\ref{theorem_LODEs_4D} LODEs in~(\ref{eq:annihilating-x}) and~(\ref{eq:annihilating-lambda}), respectively. 

We conclude the proof of Theorem~\ref{theorem_LODEs_4D} by outlining the derivation of the Pfaffian system in~(\ref{eq:pfaffian}) and~(\ref{eq:pfaffian_wrt_lambda}). On the one hand, the operator $ l_6 $ from~(\ref{basis-lambda}) readily yields~(\ref{eq:pfaffian_wrt_lambda}). On the other hand, using a multivariate analogy to transforming a higher-order LODE into a system of first order LODEs yields from operators $ l_3 $--$l_7 $ the system in~(\ref{eq:pfaffian}) --- see\cite[Sec.~6.2,~6.3]{hibi_book_13}  for an example. 

\section{Proofs Required in Section~\ref{sec:global-hgm_all}}
\subsection{Proof of Theorem~\ref{th:not stabile-lambda-2017-12-05}}
\label{section_proof_theorem_not stabile-lambda-2017-12-05}
\label{proof:th:stabile-LODE-2017-12-05}

It follows from classical results on LODE local solutions (see\cite{barkatou2009} and its references) that
the solution space on a sector of $x=\infty$ is spanned by 
functions like
$ e^{P(y)} y^q Q(\log y) \left[1 + R(y) \right]$,
where\footnote{Here, again, the generic variable symbol $ x $ stands for either of the variables of $ H^k_n(x,\lambda) $.} $y=x^{1/p}$ with $ p $ some positive integer, $P$ and $ Q$ are polynomials, and $ R $ is a holomorphic function\footnote{Differentiable in a neighborhood of any point of the domain.}.
Algorithms constructing $p, q, P, Q, R $ have been studied (see\cite{hoeij1997} and its references) and have been implemented
in computer algebra systems.
We have used {\tt formal\_sol} from Maple's DEtools to construct local solutions
for~(\ref{eq:annihilating-lambda}). Thus, we have found that the following functions span the solution space
when $\lambda \rightarrow +\infty$:
\begin{eqnarray}
h_1(x, \lambda)&=&(x\lambda)^{-1/2(1/2+n)} \exp(-2 (x\lambda)^{1/2}) (1 + {\mathcal{O}}(1/\lambda^{1/2})),\\
h_2(\lambda)&=&\lambda^{-k-1}(1+{\mathcal{O}}(1/\lambda)), \\
h_3(x, \lambda)&=&(x\lambda)^{-1/2(1/2+n)} \exp(2 (x\lambda)^{1/2}) (1 + {\mathcal{O}}(1/\lambda^{1/2})), \\
h_4(\lambda)&=&\lambda^{1-n+k} e^{\lambda} (1+{\mathcal{O}}(1/\lambda)),
\end{eqnarray}
with $h_4$ being the dominant solution.

Now, using~(\ref{equation_Hnk_definition_integral}) and~(\ref{equation_0F1_approximation}) it can readily be shown that, given $ x $, for $ \lambda \rightarrow \infty$ we have $H_n^k(x,\lambda) = {\mathcal{O}} (\exp(2 \sqrt{x \lambda}))$, i.e., there exists a constant\footnote{I.e., $ M $ is independent of $ \lambda $, but can depend on $ n $, $ k $, and $ x $.} $ M $ so that  $ | H_n^k(x,\lambda) | \le M \exp(2 \sqrt{x \lambda} ) $.
This implies that, when $ \lambda \rightarrow \infty$,
$ h_4(\lambda) \gg H_n^k(x,\lambda) $
and, hence, LODE~(\ref{eq:annihilating-lambda}) is not stabile for $H_n^k(x,\lambda)$.
This proves the first part of Theorem~\ref{th:not stabile-lambda-2017-12-05}.

On the other hand, whether or not a LODE for $ H_n^k(x,\lambda) $ is stabile when $x \rightarrow +\infty$ can be ascertained with the following Lemma, which can be proven by applying the standard procedure of the saddle-point method to~(\ref{equation_Hnk_definition_integral}). This Lemma helps prove the second part of Theorem~\ref{th:not stabile-lambda-2017-12-05}.
\begin{lemma} \label{lemma:x-asymptotic}
Given $\lambda$, when $x \rightarrow +\infty$ we have
\begin{eqnarray} \label{eq:approxHH}
  H_n^k(x,\lambda) & \sim & 2
  \frac{\Gamma(2n-1)}{\Gamma(n-1/2)} 
  (4 \sqrt{\lambda})^{-n+1/2}
  e^{\lambda}
   \exp(-P(s_0)) \left(\frac{2\pi}{P''(s_0)}\right)^{1/2}, \quad \text{with}\\
   \label{eq:s0}
 s_0 & = & \frac{1}{2}(\sqrt{\lambda} + \sqrt{\lambda+2(1+2k-n+1/2)}), \\
 \label{eq:Ps0}
 P(s) & = & (s-\sqrt{\lambda})^2-(1+2k-n+1/2)\log s.
\end{eqnarray}
\end{lemma}

Note that\footnote{We assume $\lambda+2(1+2k-n+1/2) >0$, which holds when $\lambda$ is large, so that $s_0>0$.} $s_0>0$ is obtained by solving $P'(s)=0$.
Note also that $H^k_n(x,\lambda) = {\mathcal{O}}(1)$ when
$x \rightarrow +\infty$, given $\lambda$. 

\subsection{Proof of Theorem~\ref{theorem: new-pfaffian}}
\label{section_proof_theorem_3d_LODE}

\subsubsection{LODEs w.r.t.~$ x$, given $ \lambda $}
\label{section_proof_theorem_3d_LODE_x_lambda}

On the one hand, differentiating~(\ref{equation_Hnk_definition_integral}) w.r.t.~$ x $ yields the following inhomogeneous differential equation:
    \begin{eqnarray}\label{eq:dxHkn}
        \frac{\pd{}}{\pd{}x} H^k_n(x,\lambda) = x^k e^{-x}{}_0F_1(;n;x \lambda).
    \end{eqnarray}    
On the other hand, it can be shown by using\cite[Eq.~(15.10.1), p.~394]{NIST_book_10} that ${}_0F_1(;n;z)$ satisfies the differential equation
    \begin{equation}
        \{\theta_z(\theta_z+n-1)-z\}\bullet {}_0F_1(;n;z)=0, \quad \text{with} \quad \theta_{z} = z\frac{\pd{}}{\pd{}z}.
    \end{equation}    
By letting $z=x \lambda $, according to~(\ref{eq:dxHkn}), and assuming $ \lambda $ to be a constant, the above becomes:
    \begin{equation}
    \label{0F1-LODE}
        \{\theta_x(\theta_x+n-1)-x \lambda \}\bullet {}_0F_1(;n;x \lambda )=0.
    \end{equation}    
Now,~(\ref{eq:dxHkn}) and~(\ref{0F1-LODE}) can be recast more compactly as a system of differential equations w.r.t.~$x$ for the $3$-dimensional function vector $\ug( x, \lambda) = \begin{matrix} ( H^k_n(x,\lambda)  &  {}_0F_1(;n;x \lambda) & \theta_x \bullet{}_0F_1(;n;x \lambda))^{\sf T}\end{matrix} $:
\begin{eqnarray}
\label{equation_Pfaffian_g_vector}
\frac{\pd{}}{\pd{}x}\ug=\left( 
    {\begin{array}{ccc}
        0 &  x^k e^{-x} & 0\\
        0 & 0 & 1/x\\
        0 & x \lambda & n-1\\
    \end{array}}
    \right) \ug. 
\end{eqnarray}

\subsubsection{LODEs w.r.t.~$ \varphi $, given $ \psi $}
\label{section_proof_theorem_3d_LODE_phi_psi}

If we change variables as in Remark~\ref{remark_change_variables}, i.e., define variable $  \varphi =\sqrt x$ and constant $  \psi=\sqrt \lambda $, then $ {}_0F_1(;n; x \lambda) $ becomes $ {}_0F_1(;n;  \varphi^2  \psi^2) $.
Then, on the one hand, the newly-defined function
\begin{equation}\label{eq:ut_appendix}
     v(  \varphi )=e^{-2   \varphi   \psi} {}_0F_1(;n;  \varphi^2  \psi^2)
\end{equation}
satisfies the following differential equation ensuing from~(\ref{0F1-LODE})
\begin{equation}\label{eq:new_ode1}
    \bigg\{ \theta_\varphi^2+ \left[ 4   \varphi   \psi+2(n-1)\right]\theta_ \varphi  +2 (2n-1)  \varphi   \psi \bigg\} \bullet  v(  \varphi )=0, 
\end{equation}
where $ \theta_{\varphi} $ is the differential operator $ \theta_{\varphi} = \varphi \frac{\pd{}}{\pd{}\varphi} $.
On the other hand, using~(\ref{eq:dxHkn}), identity $ \pd{} x = 2 \varphi \pd{} \varphi $, and $ {}_0F_1(;n;  \varphi^2  \psi^2) = e^{2   \varphi   \psi} v(  \varphi )  $
reveals that $ H^k_n(  \varphi ^2, \psi^2) $ satisfies the following differential equation:
\begin{equation}\label{eq:new_ode2}
    \frac{\pd{}}{\pd{}  \varphi }H^k_n(  \varphi ^2, \psi^2)=2   \varphi ^{2k+1}e^{-  \varphi ^2+2   \varphi   \psi} v(  \varphi ).
\end{equation}
Then,~(\ref{equation_LODE_system_small}) follows as the system of differential equations formed by (\ref{eq:new_ode1}), $ \theta_ \varphi $ applied to $ v(  \varphi ) $, and (\ref{eq:new_ode2}), thus proving Theorem~\ref{theorem: new-pfaffian}.

\subsection{Sketch of Proof of Theorem~\ref{th:stabile-LODE-2017-12-05}}
\label{proof_th:stabile-LODE-2017-12-05}
Given $\psi$, for $\varphi \rightarrow +\infty$ it can readily be shown that the maximum value is zero for the real part of the eigenvalues of the coefficient matrix
in~(\ref{equation_LODE_system_small_R}).
It follows from LODE local theory\cite{barkatou2009} that the dominant solution of~(\ref{equation_LODE_system_small_R}) is 
$ \mathcal{O}(1)$ w.r.t. $\varphi$.
On the other hand, Lemma~\ref{lemma:x-asymptotic} on page~\pageref{lemma:x-asymptotic} reveals that
$|| \begin{matrix} (H^k_n(  \varphi ^2, \psi^2) &  v(  \varphi ) & \theta_ \varphi  \bullet  v(  \varphi )) \end{matrix} ||$
bounds the dominant solution, completing the proof.

{\footnotesize
\bibliographystyle{IEEEtran}
\bibliography{nc_wishart,IEEEabrv,books,journals,conferences,theses,miscellaneous}}

\end{document}